\documentclass [a4paper,twoside,12pt]{article}

\usepackage{amsmath,amsfonts,amssymb}
\usepackage{vmargin,graphicx,theorem}
\usepackage{enumerate}
\usepackage[french,english]{babel}
\setpapersize[portrait]{A4}
\setmarginsrb{1cm}{2.2cm}{1cm}{1.7cm}{0cm}{0.1cm}{0.5cm}{2cm}

\newcommand{\disp}{\displaystyle}

\newcommand{\dR}{\ensuremath{\mathbb{R}}}
\newcommand{\rr}{\ensuremath{\mathbb{R}}}
\newcommand{\ee}{\ensuremath{\mathbb{E}}}

\newtheorem{ethm}{}
\newtheorem{theo}[ethm]{Theorem}

\newtheorem{ecor}[ethm]{Corollary}
\newtheorem{coro}[ethm]{Corollary}

\newtheorem{eprop}[ethm]{Proposition}
\newtheorem{prop}[ethm]{Proposition}

\newtheorem{elem}[ethm]{Lemma}
\newtheorem{lemma}[ethm]{Lemma}

\newtheorem{defi}[ethm]{Definition}

\newtheorem{erem}[ethm]{Remark}

\newcommand{\proofend}{~$\rhd$}
\newcommand{\proofbegin}{~$\lhd$}

\newenvironment{eproof}
               {\noindent {\emph{\textbf{Proof}}}\\\proofbegin~}
               {\proofend\\}

\newcommand{\p}[4]{{#3}\!\left#1{#4}\right#2}

\newcommand{\ABS}[1]{\ensuremath{{\left| #1 \right|}}} 
\newcommand{\PAR}[1]{\ensuremath{{\left(#1\right)}}} 
\newcommand{\SBRA}[1]{\ensuremath{{\left[#1\right]}}} 
\newcommand{\BRA}[1]{\ensuremath{{\left\{#1\right\}}}} 
\newcommand{\NRM}[1]{\ensuremath{{\left\Vert #1\right\Vert}}} 

\renewcommand{\phi}{\varphi}

\newcommand{\varf}[1]{\mathbf{Var}_{#1}}
\newcommand{\entf}[1]{\mathbf{Ent}_{#1}}

\newcommand{\ent}[2]{\p(){\entf{#1}}{#2}}

\newcommand{\var}[2]{\p(){\varf{#1}}{#2}}

\newcommand{\PT}[1]{\mathbf{P}_{\!#1}}
\newcommand{\Pt}[1][t]{\ensuremath{\mathbf{P_{\!#1}}}}

\begin{document}

\title{\sl Phi-entropy inequalities for diffusion semigroups }
\author{ Fran\c cois Bolley and Ivan Gentil}
\maketitle

\begin{center}
\small{Universit\'e Paris-Dauphine,
Ceremade, UMR CNRS 7534\\
Place du Mar\'echal de Lattre de Tassigny,
F-75116 Paris cedex 16\\
bolley, gentil@ceremade.dauphine.fr}
\end{center}

\bigskip

\begin{center}
{\bf Abstract}
\end{center}

We obtain and study new $\Phi$-entropy inequalities for diffusion semigroups, with Poincar\'e or logarithmic Sobolev inequalities as particular cases. From this study we derive the asymptotic behaviour of a large class of linear Fokker-Plank type equations under simple conditions, widely extending previous results. Nonlinear diffusion equations are also studied by means of these inequalities. The $\Gamma_2$ criterion of D.~Bakry and M. Emery appears as a main tool in the analysis, in local or integral forms.

\begin{center}
{\bf R\'esum\'e}
\end{center}
Nous obtenons et \'etudions une nouvelle famille d'in\'egalit\'es $\Phi$-entropiques pour des semigroupes de diffusion, incluant les in\'egalit\'es de Poincar\'e et de Sobolev logarithmiques. Nous en d\'eduisons le comportement en temps grand des solutions d'une grande classe d'\'equations lin\'eaires de type Fokker-Planck, sous de simples conditions. Nous \'etudions \'egalement certaines  \'equations de diffusion nonlin\'eaires \`a l'aide de ces in\'egalit\'es. Cette \'etude utilise de mani\`ere cruciale le crit\`ere  $\Gamma_2$ de D.~Bakry et M.~Emery, sous des formes locales et int\'egrales.


\bigskip

\bigskip


\bigskip

\section*{Introduction}
Functional inequalities such as the Poincar\'e inequality and the logarithmic Sobolev inequality of L.~Gross  have revealed adapted to obtain estimates in  the asymptotic in time behaviour of diffusion Markov semigroups and of solutions to Fokker-Planck type equations for instance. These two inequalities respectively imply an exponential decay in time of the (relative) variance and of the Boltzmann entropy of the solution. A natural interpolation between the Poincar\'e and  logarithmic Sobolev inequalities, via an interpolation between the variance and the entropy, consists in the generalized Poincar\'e or Beckner inequalities, which translate into an exponential decay of the $L^p$-norms of the solutions for $p$ between 1 and~2. These interpolating inequalities are part of the family of $\Phi$-entropy inequalities, where $\Phi$ belongs to a class of convex functions satisfying additional admissibility assumptions:
 they correspond to the instances of maps $\Phi$ given by $\Phi(x) = x^p$ with $p$ between $1$ and $2.$

Derivations of such functional inequalities are to a large extent based on the so-called $\Gamma_2$ criterion introduced by D. Bakry and M. Emery. This criterion is a local condition on the coefficients of the infinitesimal generator of the semigoup, or  of the Fokker-Planck equation. It is a sufficient condition to $\Phi$-entropy inequalities for the possible ergodic measure of the semigroup, which gives its long time behaviour;  but it is a necessary and sufficient condition to such  inequalities for the associated Markov kernel of the semigroup (at each time $t$). On the other hand, for the Poincar\'e and logarithmic Sobolev inequalities, it can be replaced by a weaker nonlocal condition, called integral criterion, in the study of the sole ergodic measure.

\bigskip

In this work we consider a general diffusion semigroup on $\rr^n$. A first section is devoted to a general and simplified presentation and derivation of  $\Phi$-entropy inequalities for general admissible $\Phi$'s, both for the Markov kernel and the possible ergodic measure of the semigroup. A general integral criterion is obtained, which extends the Poincar\'e and logarithmic Sobolev inequality cases, and interpolates between them when specified to the Beckner power law cases  $\Phi(x) = x^p$ with $1 < p <2$ (see Proposition~\ref{prop-ci}). We finally study the asymptotic in time behaviour of a large class of diffusion semigroups: rephrased in the Fokker-Planck type equation setting, we show how to simply obtain the existence of a unique stationary state and the convergence of all the solutions towards it, in $\Phi$-entropy senses, and with a precise rate (see Theorem~\ref{theo-cvfp}). The method applies to a much wider class of linear equations than those previoulsy studied, in which the arguments are strongly based on a deeper knowledge of the limit measure, and in particular on its explicit expression.

In Section 2 we focus on power law entropies $\Phi(x) = x^p$ with $1< p <2$: we strenghten the Beckner inequalities by deriving and studying certain power law $\Phi$-entropy inequalities, introduced by A.~Arnold and J.~Dolbeault, both for the Markov kernel and the possible ergodic measure (see Theorem~\ref{thm-main}). Inequalities for the Markov kernel are shown to be equivalent to the local $\Gamma_2$ criterion, and the corresponding inequalities for the ergodic measure are implied by a weaker and adapted integral condition, thus improving on the results by A.~Arnold and J.~Dolbeault. We study properties of these inequalities, proving that they constitute a new monotone interpolation between the Poincar\'e and the logarithmic Sobolev inequalities (see Proposition~\ref{prop-7} and Remark~\ref{lsibeckner+}).

In Section 3 we show that the $\Phi$-entropy inequalities may not hold for nonadmissible functions~$\Phi$; however, for $\Phi(x)=x^p$ with $p$~positive or $\Phi(x) = x \ln x$ we obtain similar inequalities but with an extra term (see Theorem~\ref{thm-5}).

As an application we show in a last section how  the local functional inequalities for Markov semigroups obtained so far can be extended to inhomogeneous semigroups (see Proposition \ref{inhom}) and imply analogous properties on solutions  of an instance of nonlinear Fokker-Planck evolution equation in a very simple way (see Theorem \ref{theo-mkv}).

 \section{Phi-entropies}\label{sectphi}

We consider a Markov semigroup $(\PT{t})_{t \geq 0}$ on $\rr^n$, acting on functions on $\rr^n$ by
$$
\PT{t} f(x) = \int_{\rr^n} f(y) \, p_t(x, dy)
$$
for $x$ in $\rr^n.$ The kernels $p_t(x, dy)$ are probability measures on $\rr^n$ for all $x$ and $t \geq 0$, called transition kernels. Moreover we assume that the Markov infinitesimal generator $L = \disp \frac{\partial}{\partial t} \Big\vert_{t=0^+} \PT{t}$ is given by
$$
L f (x)= \sum_{i,j=1}^n D_{ij} (x) \frac{\partial^2 f}{\partial x_i \partial x_j} (x) - \sum_{i=1}^n a_i (x)  \frac{\partial f}{\partial x_i} (x)
$$
where  $D(x) = (D_{ij}(x))_{1 \leq i,j \leq n}$ is a symmetric $n \times n$ matrix, nonnegative in the sense of quadratic forms on $\rr^n$ and with smooth coefficients; also $a(x) = (a_i(x))_{1 \leq i \leq n}$ has smooth coefficients. Such a semigroup or generator  is called a {\it diffusion}, and we refer to  \cite{bakrystflour}, \cite{ledouxmarkov}  or  \cite {bakrytata} for backgrounds on these semigroups and forthcoming notions. 

\bigskip

If $\mu$ is a Borel probability measure on $\rr^n$ and $f$ a $\mu$-integrable Borelian function on $\rr^n$ we let
$$
\mu(f) = \int_{\rr^n} f(x) \, \mu(dx).
$$
If moreover $\Phi$ is a convex function on an interval $I$ of $\rr$ and $f$ an $I$-valued Borelian  function such that $f$ and $\Phi(f)$ be $\mu$-integrable,  we let 
$$
\entf{\mu}^\Phi\PAR{f} = \mu(\Phi(f)) - \Phi(\mu(f)) 
$$
be the $\Phi$-entropy of $f$ under $\mu$ (see  \cite{chafai04} for instance).
Two fundamental examples are  $\Phi(x) = x^2$ on $\rr$, for which we let $\var{\mu}{f} = \entf{\mu}^\Phi\PAR{f}$ be the variance of $f,$ and $\Phi(x) = x \ln x $ on $]0, +\infty[$, for which we let $\ent{\mu}{f} = \entf{\mu}^\Phi\PAR{f}$ be the Boltzmann entropy of a positive function $f.$
 
 By Jensen's inequality  the $\Phi$-entropy $\entf{\mu}^\Phi (f)$ is always a nonnegative quantity, and in this work we are interested in deriving upper and lower bounds on  $\entf{\mu}^\Phi\PAR{f}$ or $\entf{\PT{t}}^\Phi\PAR{f} (x)  = \entf{p_t(x,dy)}^\Phi\PAR{f} = \Pt \Phi(f) (x) - \Phi(\Pt f) (x)$ under adequate assumptions on $L$. Applications of such bounds  to inhomogeneous Markov semigroups and to linear and nonlinear Fokker-Planck type equations will also be given, in one of which the diffusion matrix $D$ being nonnegative but not (strictly) positive.
 
 \bigskip
 
 \subsection{The carr\'e du champ and $\Gamma_2$ operators}
 Bounds on $\entf{\PT{t}}^\Phi\PAR{f}$ and assumptions on $L$ will be given in terms of the {\it carr\'e du champ} operator associated to $L$, defined by
 $$
 \Gamma(f,g) = \frac{1}{2} \Big( L(fg) - f \, Lg - g \, Lf \Big).
 $$
For simplicity we shall let $\Gamma(f) = \Gamma(f,f).$ 
 Assumptions on $L$ will also be given in terms of the $\Gamma_2$ operator defined by
 $$
  \Gamma_2(f) = \frac{1}{2} \Big( L \Gamma(f) - 2\Gamma (f, Lf) \Big).
 $$
 
 \begin{defi}
  If $\rho$ is a real number, we say that the semigroup $(\PT{t})_{t \geq 0}$ (or the infinitesimal generator $L$) satisfies the {\it $CD(\rho,\infty)$  criterion} if
$$
\Gamma_2(f)\geq \rho \, \Gamma(f)
$$
for all functions $f$.
\end{defi}

This criterion is a special case of the curvature-dimension criterion $CD(\rho,m)$ with $\rho\in\dR$ and $m\geq 1$ proposed by D. Bakry and M. Emery (see \cite{bakryemery}).
\bigskip

\noindent {\bf Example 1} A fundamental example is the heat semigroup on $\dR^n$ defined by
$$
\PT{t}f (x) = \int_{\rr^n} f(y) \frac{e^{-\frac{\NRM{x-y}^2}{4t}}}{\PAR{4\pi t}^{n/2}}\, dy.
$$
Its generator is the Laplacian and it satisfies the $CD(0,\infty)$ criterion.

\bigskip

\noindent {\bf Example 2}
Another fundamental  example is the Ornstein-Uhlenbeck semigroup defined by
$$
\PT{t}f (x) = \int_{\rr^n} f(e^{-t} x + \sqrt{1 - e^{-2t}} y) \, \gamma(dy)
$$
where $\gamma(dy) = (2\pi)^{-n/2} \exp(-{\NRM{y}^2}/{2}) \, dy$ is the standard Gaussian measure on $\rr^n.$ Its infinitesimal generator is given by
$$
Lf(x) = \Delta f(x) - < x , \nabla f(x) >
$$
where $\Delta, \nabla$ and $< \cdot , \cdot>$ respectively stand for the Laplacian and gradient operators and the scalar product on $\rr^n.$ Then  the carr\'e du champ and $\Gamma_2$ operators are  given by
$$
\Gamma (f) = \Vert \nabla f \Vert^2
$$
and
$$
\Gamma_2 (f) = \Vert \textrm{Hess} f \Vert_2^2 + \Vert \nabla f \Vert^2
$$
where Hess$f$ is the Hessian matrix of $f$ and $\Vert M \Vert_2^2 = \displaystyle \sum_{i,j=1}^n M_{ij}^2$ if $M$ is the matrix $(M_{ij})_{1 \leq i,j \leq n}.$ In particular the Ornstein-Uhlenbeck semigroup satisfies the $CD(1,\infty)$ criterion. 

\bigskip

 The carr\'e du champ  associated to a general infinitesimal generator
$$
Lf (x) = \sum_{i,j=1}^n D_{ij} (x) \frac{\partial^2 f}{\partial x_i \partial x_j} (x) - \sum_{i=1}^n a_i (x) \frac{\partial f}{\partial x_i} (x)
$$
is given by
$$
\Gamma(f) (x) = < \nabla f (x), D(x) \,  \nabla f(x) >.
$$
Expressing $\Gamma_2$ is more complex: For instance,

$\bullet$  if $D$ is constant then
$$
\Gamma_2(f) (x) = \textrm{trace} \Big((D \, \textrm{Hess} f (x)) ^2 \Big) +  < \nabla f (x), J a (x) \, D \,  \nabla f(x) >
$$
where $\disp J a = \Big( \frac{\partial a_i}{\partial x_j} \Big)_{i,j}$ is the {\it Jacobian matrix} of $a$; then $L$ satisfies the $CD(\rho,\infty)$  criterion if and only if
\begin{equation}\label{cdrhocst}
\frac{1}{2} \big( J a (x) D + (J a (x) D )^* \big) \geq \rho \, D
\end{equation}
for all $x$ as quadratic forms on $\rr^n$, where $M^*$ is the transposed matrix of a matrix $M;$

$\bullet$ if $D(x) = d(x) I$ is a scalar matrix then, letting $\partial_i g = \disp  \frac{\partial g}{\partial x_i}$ and $\partial^2_{ij} g = \disp \frac{\partial^2 g}{\partial x_i \partial x_j}$,
$$
\Gamma_2(f)  = \sum_{i =1}^n  [d \, \partial^2_{ii} f +  \partial_i d \, \partial_i f - \frac{1}{2} \sum_{k=1}^n \partial_k d \, \partial_k f]^2  
+  \sum_{i \neq j} [d \,  \partial^2_{ij} f + \frac{1}{2} (\partial_i d \, \partial_j f + \partial_j d \, \partial_i f)]^2 + \sum_{i,j= 1}^n \partial_i f M_{ij} \partial_j f
$$
where
$$
M = \frac{1}{2} (d \, \Delta d -  < a , \nabla d  >-\NRM{\nabla d}^2) I + (\frac{1}{2} - \frac{n}{4}) \nabla d \otimes \nabla d + d^2 \, J a \, ;
$$
then $L$, which is $d(x) \Delta - <a, \nabla >$  in this case, satisfies the $CD(\rho,\infty)$  criterion if and only if 
\begin{equation}\label{cdrhoscalaire}
\frac{1}{2} \big( M(x) + M(x)^* \big)  \geq \rho \, d(x) \, I
\end{equation}
for all $x$, as quadratic forms on $\rr^n$. This condition can also be found in \cite{amtucpde01}, as will be discussed more in detail in Section~\ref{subsectFP}.

\smallskip

The $CD(\rho, \infty)$ criterion for a general $L$ with {\it positive} diffusion matrix $D$ is discussed in \cite{arnoldcarlenju08}.

\bigskip
\subsection{Poincar\'e and logarithmic Sobolev inequalities}\label{subsectPLSI}

The $CD(\rho,\infty)$ criterion for a $\rho\in\dR$ is well adapted to deriving upper and lower bounds on $\Phi$-entropies.

\noindent
{\bf Example 3 (\cite{bakrystflour})} For $\Phi(x) = x^2$, the following four assertions are equi\-valent, with $\disp  \frac{1 - e^{-2\rho t}}{\rho}$ and $\disp \frac{e^{2\rho t} -1}{\rho}$ replaced by $2 t$ if $\rho =0$:

\begin{enumerate}[(i)]
\item
 the semigroup $(\PT{t})_{t \geq 0}$ satisfies the ${CD}(\rho,\infty)$ criterion;
\item
the semigroup $(\PT{t})_{t \geq 0}$ satisfies the commutation relation
$$
\Gamma (\Pt f) \leq e^{-2\rho t} \, \Pt (\Gamma (f))
$$
for all positive $t$ and all functions $f$;
\item 
the semigroup $(\PT{t})_{t \geq 0}$ satisfies the local Poincar\'e inequality
 \begin{equation}\label{localPI}
 \var{\Pt}{f} \leq \frac{1 - e^{-2\rho t}}{\rho} \, \Pt(\Gamma (f))
 \end{equation}
  for all positive $t$ and all functions $f$; 
 \item
 the semigroup $(\PT{t})_{t \geq 0}$ satisfies the reverse local Poincar\'e inequality
  \begin{equation}\label{reverselocalPI}
 \var{\Pt}{f} \geq \frac{e^{2\rho t} -1}{\rho} \, \Gamma (\Pt f)
  \end{equation}
  for all positive $t$ and all functions $f$.
\end{enumerate}

Let us note that  (iii) and (iv) together imply (ii), and that in fact (i) for instance holds as soon as (iii) or (iv) holds for all $t$ in a neighbourhood of $0.$

\medskip

\noindent
{\bf Example 4 (\cite{bakrystflour})} For $\Phi(x) = x \ln x$, the following four assertions are equivalent, with the same convention for $\rho =0$:

\begin{enumerate}[(i)]
\item
 the semigroup $(\PT{t})_{t \geq 0}$ satisfies the ${CD}(\rho,\infty)$ criterion;
\item
the semigroup $(\PT{t})_{t \geq 0}$ satisfies the commutation relation
$$
\Gamma(\Pt f) \leq e^{-2\rho t}{\Pt \Big(\sqrt{\Gamma(f)} \Big)}^2
$$
for all positive $t$ and all positive functions $f$; 
\item 
the semigroup $(\PT{t})_{t \geq 0}$ satisfies the local logarithmic Sobolev inequality
\begin{equation}\label{localLSI}
\ent{\Pt}{f} \leq \frac{1 - e^{-2\rho t}}{2\rho} \, \Pt\PAR{\frac{\Gamma (f)}{f}}
\end{equation}
for all  positive $t$ and all positive functions $f$; 
 \item
the semigroup $(\PT{t})_{t \geq 0}$ satisfies the reverse local logarithmic Sobolev inequality
$$
 \ent{\Pt}{f} \geq \frac{ e^{2\rho t}-1}{2\rho} \, {\frac{\Gamma (\Pt f)}{\Pt f}}
$$ 
for all  positive $t$ and all positive functions $f$.
\end{enumerate}

Let us note that, contrary to this second case, the equivalences in Example 3 when $\Phi(x) = x^2$ hold in a more general setting, when the generator $L$ is not a diffusion semigroup.

\bigskip

A Borel probability measure $\mu$ on $\rr^n$ is called {\it invariant} for the semigroup $(\Pt)_{t \geq 0}$ if $\mu(\Pt f) = \mu(f)$ for all $t$ and $f$, or equivalently if $\mu(L f)=0$ for all $f.$ Then we say that the semigroup $(\Pt)_{t \geq 0}$ is {\it $\mu$-ergodic} if $\Pt f$ converges to
 $\mu(f)$ as $t$ tends to infinity,  in $L^2(\mu)$ for all functions $f$. 
  For instance, if $\mu$ is an  invariant probability measure for the semigroup $(\Pt)_{t \geq 0},$ then $(\Pt)_{t \geq 0}$ is $\mu$-ergodic as soon as the carr\'e du champ $\Gamma$ vanishes only on constant functions, that is, for such diffusion semigroups, as soon as the matrix $D(x)$ is positive for all $x$.
 
Then let  $\mu$ be an ergodic probability measure for $(\Pt)_{t \geq 0}$;  if the $CD(\rho, \infty)$ criterion holds with $\rho >0$, it follows from \eqref{localPI} and \eqref{localLSI} respectively that the measure $\mu$ satisfies the Poincar\'e inequality
 \begin{equation}\label{PI}
 \var{\mu}{f} \leq \frac{1}{\rho} \mu(\Gamma(f))
 \end{equation}
 for all functions $f$
 and the logarithmic Sobolev inequality
\begin{equation}\label{LSI}
\ent{\mu}{f} \leq  \frac{1}{2 \rho} \mu\PAR{\frac{\Gamma(f)}{f}}
\end{equation}
 for all positive functions $f$.
The logarithmic Sobolev inequality \eqref{LSI} with constant $1/2 \rho$, introduced by L. Gross in \cite{gross75} (see also~\cite{gross93}), is known to imply the Poincar\'e inequality \eqref{PI} with constant $1/\rho.$

 In fact \eqref{PI} and \eqref{LSI} hold under a condition on the generator $L$ weaker than the $CD(\rho, \infty)$ criterion (see \cite{bakryemery,ledoux92,logsob} for instance): First of all, if  $(\Pt)_{t \geq 0}$ is $\mu$-ergodic and $\rho$ is a positive number then $\mu$ satisfies the Poincar\'e inequality \eqref{PI} if it satisfies the averaged $CD(\rho, \infty)$ condition called {\it integral criterion}
 \begin{equation}
 \label{critereint}
\mu(\Gamma_2(f)) \geq \rho \, \mu(\Gamma(f))
\end{equation}
for all $f$.  If moreover $\mu$ is {\it reversible} with respect to the semigroup $(\Pt)_{t \geq 0}$, that is, if  $\mu(f \Pt g) = \mu(g \Pt f)$ for all functions  $f$ and $g$, then the integral criterion~\eqref{critereint} is {\it equivalent} to the Poincar\'e inequality~\eqref{PI}.  

On the other hand,  the logarithmic Sobolev inequality \eqref{LSI} is  implied by  the integral criterion
\begin{equation}\label{criteresuperint}
\mu(e^f \, \Gamma_2(f)) \geq \rho \, \mu(e^f \, \Gamma(f))
\end{equation}
for all $f.$ This result is often stated under the reversibility condition, which is useless for diffusion semigroups, as we shall see in Proposition~\ref{prop-ci}. As pointed out by B. Helffer (see \cite[p.~114]{helffer02} or \cite[p.~91]{logsob}), the converse does not hold, even under the reversibility condition.

\medskip

\bigskip
\subsection{$\Phi$-entropy inequalities}
\label{sub-11}
(Local) Poincar\'e and logarithmic Sobolev inequalities for the semigroup $(\Pt)_{t \geq 0}$ are part of a large family of functional inequalities introduced in \cite{bakryemery} and developed in \cite{chafai04} and \cite{bakrytata}:

\begin{theo}\label{theophisobolev}
Let $\rho$ be a real number and $\Phi$ be a $C^4$ strictly convex function on an  interval $I$ of $\rr$ such that $-1/\Phi''$ be convex. Then the following three assertions are equivalent, with $\disp  \frac{1 - e^{-2\rho t}}{2 \rho}$ and $\disp \frac{e^{2\rho t} -1}{2 \rho}$ replaced by $t$ if $\rho =0$:
\begin{enumerate}[(i)]
\item
 the semigroup $(\PT{t})_{t \geq 0}$ satisfies the ${CD}(\rho,\infty)$ criterion;
\item 
 the semigroup $(\PT{t})_{t \geq 0}$ satisfies the local $\Phi$-entropy inequality
 \begin{equation}\label{localPHII}
 \entf{\Pt}^\Phi (f) \leq \frac{1 - e^{-2\rho t}}{2 \, \rho} \,  \Pt( \Phi''(f) \Gamma (f))
 \end{equation}
   for all positive $t$ and all $I$-valued functions $f$; 
 \item
the semigroup $(\PT{t})_{t \geq 0}$ satisfies the reverse local $\Phi$-entropy inequality
  \begin{equation}\label{reverselocalPHII}
 \entf{\Pt}^\Phi (f) \geq \frac{e^{2\rho t} -1}{2 \, \rho} \, \Phi''(\Pt f) \, \Gamma (\Pt f)
  \end{equation}
    for all positive $t$ and all $I$-valued functions $f$.
 \end{enumerate}
 
If moreover the probability measure $\mu$ is ergodic for the semigroup $(\Pt)_{t \geq 0}$, and $\rho>0$,  then $\mu$ satisfies the $\Phi$-entropy inequality 
\begin{equation}\label{PHII}
 \entf{\mu}^\Phi (f)   \leq \frac{1}{2 \rho} \, \mu (\Phi''(f) \, \Gamma (f))
\end{equation}
for all $I$-valued functions $f$.
\end{theo}

For future use, we shall give a proof of Theorem~\ref{theophisobolev} slightly different from the one given in \cite{chafai04} and  \cite{bakrytata}. Before doing so we make some comments on this result.

\bigskip

A function $\Phi$ satisfying the conditions of Theorem~\ref{theophisobolev} will be called an {\it admissible function}.  The functions $\Phi: x \mapsto x^2$ or more generally $a x^2 + b x + c$ on $\rr$ and $x \mapsto x \ln x$ or more generally $(x+a) \ln (x+a) + b x +c$ on $]-a, +\infty[$ are the solutions to $(1/\Phi'')''=0$ and thus are admissible. They respectively lead to the (local) Poincar\'e and logarithmic Sobolev inequalities of Section~\ref{subsectPLSI}. More generally, for any $1 \leq p \leq 2$ the function
\begin{equation}
\label{eq-phip}
\Phi_p(x)= \left\{
\begin{array}{cl}
\disp \frac{x^p- x}{p(p-1)}, \quad x > 0 & \text{ if }p\in ]1,2]\\
\disp x \ln x, \quad x>0 & \text{ if }p=1
\end{array}
\right.
\end{equation}
is admissible. For this entropy $\Phi_p$ with $p$ in $]1,2]$ the $\Phi$-entropy inequality \eqref{PHII}  writes
\begin{equation}\label{beckneravecp}
\frac{\mu(f^p) - \mu(f)^p}{p(p-1)} \leq \frac{1}{2 \rho} \, \mu (f^{p-2} \, \Gamma (f))
 \end{equation}
for all positive functions $f$, or equivalently
\begin{equation}\label{beckner}
\frac{\mu(g^2) - \mu(g^{2/p})^p}{p-1} \leq \frac{2}{ p  \rho} \, \mu (\Gamma (g))
 \end{equation}
for all positive functions $g$, with $g = f^{p/2}.$
These $\Phi$-entropy inequalities have been studied in \cite{beckner89} for the uniform measure on the sphere and the Gaussian measure, and are called {\it generalized Poincar\'e} or {\it Beckner's inequalities}. For given $g>0$ the map $p \mapsto \mu(g^{2/p})^p$ is convex on $]0, +\infty[$, so that the quotient $\disp \frac{\mu(g^2) - \mu(g^{2/p})^p}{p-1}$ is nonincreasing with respect to $p$, $p\neq1$ (see~\cite{latalaoles}). Moreover its limit for $p$ tending to $1$ is $\ent{\mu}{g^2}$, so that
$$
\disp \frac{\mu(g^2) - \mu(g^{2/p})^p}{p-1} \leq \var{\mu}{g} \leq \disp \frac{\mu(g^2) - \mu(g^{2/q})^q}{q-1}
\leq \ent{\mu}{g^2} \leq \disp \frac{\mu(g^2) - \mu(g^{2/r})^r}{r-1}
$$
for all positive functions  $g$ and all $p,q,r$ such that  in $0<r<1<q<2<p$: in this sense the Beckner inequalities~\eqref{beckner} for $p$ in $]1,2]$ give a natural interpolation between the weaker Poincar\'e inequality~\eqref{PI} for positive functions, and then for all functions, and the stronger logarithmic Sobolev inequality~\eqref{LSI}, with~\eqref{beckner}  being \eqref{PI} for $p=2$ and giving \eqref{LSI} in the limit $p \to 1$ (see \cite{ledouxmarkov} for instance).

\begin{erem}\label{remphiadmi}
A general study of admissible functions is performed in \cite{amtucpde01}, \cite{chafai04, chafai06} and \cite{evol}: for instance a $C^4$ strictly convex function $\Phi$ on an interval $I$ of $\rr$ is admissible if and only if $\Phi^{(4)}(x) \, \Phi'' (x) \geq 2 \, \Phi'''(x)$ for all $x$, and if and only if the map $(x,y) \mapsto \Phi''(x) \, y^2$ is convex on $I \times \rr.$

In fact one can note that a $C^4$ function $\Phi$ on $I$ is admissible if and only if $1/\Phi''$ is a $C^2$ positive concave function on $I.$ 

First of all, this leads to other examples of admissible functions, such as:

\begin{itemize}
\item 
 if  $\alpha\in[1,2[$ and  $\beta\in\dR$ then one can find $a\geq 1$ such that the map $\Phi$ defined by 
$$
\Phi(x)=(x+a)^\alpha\PAR{\ln(x+a)}^\beta
$$
 is admissible on  $[0,+\infty[$, thus extending the family of Beckner's entropies;
\item if $a$ is a positive real number, then a primitive of the function $x\mapsto  \ln (e^{ax}-1)$ is also an admissible function on $]0, +\infty[.$
\end{itemize}

Then it enables to recover the fact that  the set of admissible functions on a given interval $I$ is a convex vector cone, as pointed out  in \cite[Remark 5]{chafai04}).

Indeed, let $\Phi_1$ and $\Phi_2$ be two admissible functions and $\lambda\in[0,1].$ Then $\lambda\Phi_1+(1-\lambda)\Phi_2$ is convex and $1/(\lambda\Phi_1''+(1-\lambda)\Phi_2'')$ is concave since
$$
\PAR{ \Big({\frac{1}{\theta_1}+\frac{1}{\theta_2} \Big)}^{\! \! -1} }''=\theta_1^6\theta_2^6\frac{-2(\theta_1\theta_2'-\theta_1'\theta_2)^2+\theta_1^3\theta_2''+\theta_1''\theta_2^3+\theta_1\theta_1''\theta_2^2+\theta_1^2\theta_2\theta_2''}{\PAR{\theta_1+ \theta_2}^3}
$$
where  $\theta_1=1/(\lambda\Phi_1'')$ and $\theta_2=1/((1-\lambda)\Phi_2'')$ are positive concave functions.  Hence $\lambda\Phi_1+(1-\lambda)\Phi_2$ is an admissible function. 
\end{erem}

\bigskip

\textbf{\textit{Proof of Theorem~\ref{theophisobolev}}}

\noindent \proofbegin\,
 We first assume $(i)$ and prove $(ii)$ and $(iii)$. We let $t >0$ be fixed and we consider the function
 \begin{equation}
\label{eq-psi}
\psi\PAR{s} =  \PT{s}\PAR{\Phi\PAR{\PT{t-s} f}}
\end{equation}
 so that
 $$
 \entf{\Pt}^\Phi(f)=\Pt(\Phi(f)) - \Phi(\Pt f) = \psi(t) - \psi(0).
 $$
 
 Let us first  admit the following
 \begin{lemma}
 \label{lemmepsi'}
 For any $C^4$ function $\Phi$ with nonvanishing second derivative the function 
 $\psi(s) = \PT{s}\PAR{\Phi\PAR{\PT{t-s} f}}$ is twice differentiable on $[0,t]$, with
$$ 
 \psi'(s)=\PT{s}\PAR{\Phi''(\PT{t-s}   f) \Gamma(\PT{t-s}   f)} = \PT{s}\PAR{\frac{\Gamma(\Phi'(\PT{t-s}   f))}{\Phi''(\PT{t-s}   f)}}
 $$
 and
\begin{equation}
\label{eq-seconde}
\psi''(s)=2 \, \PT{s}\PAR{\frac{\Gamma_2(\Phi'(\PT{t-s}   f))}{\Phi''(\PT{t-s}   f)}}\\
+ \PT{s}\PAR{\PAR{\frac{\Gamma (\Phi'(\PT{t-s}   f))}{\Phi''(\PT{t-s}   f)}}^2\PAR{\frac{-1}{\Phi''}}''(\PT{t-s}f)}.
\end{equation}
\end{lemma}

\medskip

By assumption on $\Phi$ the second term on the right hand side of \eqref{eq-seconde} is nonnegative, so 
$$
\psi''(s)\geq 2\rho\psi'(s), \quad s \in [0,t]
$$
by the $CD(\rho,\infty)$ criterion.  Now for $0 \leq u \leq v \leq t$ we integrate over  $[u,v]$ to obtain 
$$
\psi'(u) \leq \psi'(v) e^{2\rho(u-v)}. 
$$

For $u=s$ and $v=t$, integrating in $s$  over the set $[0,t]$ yields
$$
\psi(t) - \psi(0) \leq \psi' (t) \, \frac{1 - e^{-2\rho t}}{\rho},
$$
which is $(ii)$. 

For $u=0$ and $v=s$, integrating in $s$ over the set $[0,t]$ yields
$$
\psi(t) - \psi(0) \geq \psi' (0) \, \frac{e^{2\rho t} -1}{\rho},
$$
which is $(iii).$

\medskip

Let us conversely assume that $(ii)$ or $(iii)$ holds. For $f = 1 + \varepsilon g$ the left hand side in \eqref{localPHII} and \eqref{reverselocalPHII} is 
$$
 \entf{\Pt}^\Phi (f) = \frac{\varepsilon^2}{2} \, \Phi''(1) \, \var{\Pt} {g} + o(\varepsilon^2)
$$
and the right hand side is given by
$$
\Pt (\Phi''(f) \Gamma (f))) = \varepsilon^2 \, \Phi''(1) \PT{t} (\Gamma (g)) + o(\varepsilon^2).
$$
Hence, in the limit $\varepsilon \to 0,$ $(ii)$ implies the local Poincar\'e inequality \eqref{localPI} and $(iii)$ implies the reverse local Poincar\'e inequality \eqref{reverselocalPI} which are equivalent to the $CD(\rho, \infty)$ criterion.
This concludes the proof of Theorem \ref{theophisobolev}. 
\proofend

\bigskip

We now turn to the \emph{\textbf{proof of Lemma~\ref{lemmepsi'}}}.

\noindent \proofbegin$\,$
The first derivative of $\psi$ is given by
$$
\psi'(s)=\PT{s}\PAR{L \Phi(\PT{t-s}   f) - \Phi'(\PT{t-s}   f) \, L \PT{t-s} f}.
$$
But $L$ is a diffusion so satisfies the identities
\begin{equation}
\label{eq-diffusion}
L\Phi(g)=\Phi'(g)Lg+\Phi''(g)\Gamma(g)\,\,\rm{and}\,\, 
\Gamma(\Phi'(g))=\Phi''^2(g)\Gamma(g)
\end{equation}
(see for instance \cite[Lemme~1]{bakryemery} or \cite[p. 31]{logsob}) Hence
$$ 
\psi'(s)=\PT{s}\PAR{\Phi''(\PT{t-s}   f) \Gamma(\PT{t-s}   f)} = \PT{s}\PAR{\frac{\Gamma(\Phi'(\PT{t-s}   f))}{\Phi''(\PT{t-s}   f)}}.
$$

Then the derivative of $\psi'$ is
$$ 
\psi''(s)=\PT{s}\BRA{L\PAR{\Phi''(g)\Gamma(g)}-\Phi'''(g)Lg \, \Gamma(g)-2\Phi''(g)\Gamma(g,Lg)}
$$
where $g=\PT{t-s}f$. Then the definition 
$$
L(f_1f_2)=2\Gamma(f_1,f_2)+f_1Lf_2+f_2Lf_1
$$
of $\Gamma$, the identities \eqref{eq-diffusion} and the definition of $\Gamma_2$ yield
$$ 
\psi''(s)=\PT{s}\BRA{\frac{1}{\Phi''(g)}\SBRA{2\Phi'''(g)\Phi''(g)\Gamma(g,\Gamma(g))+2\Phi''^2(g)\Gamma_2(g)+\Phi^{(4)}(g)\Phi''(g)\Gamma(g)^2}}. 
$$
But again $L$ is a diffusion operator, so satisfies the identity 
$$
\Gamma_2(\Phi'(g))=(\Phi''(g))^2\Gamma_2(g)+\Phi''(g)\Phi'''(g)\Gamma(g,\Gamma(g))+\Phi'''^2(g)\Gamma(g)^2,
$$ 
for all functions $g$ (see for example \cite[Lemme~3]{bakryemery} or \cite[Lemme~5.1.3]{logsob}), which gives the expression of the second derivative of $\psi$.
\proofend

\medskip
\bigskip

As for the Poincar\'e and logarithmic Sobolev inequalities of Section~\ref{subsectPLSI}, the pointwise $CD(\rho,\infty)$ criterion can be replaced by an integral criterion to get the $\Phi$-entropy inequality~\eqref{PHII}:
\begin{eprop}
\label{prop-ci}
Let $\rho$ be a positive number and $\Phi$ be an admissible function on an interval $I$. If the probability measure $\mu$ is ergodic for the diffusion semigroup $(\Pt)_{t \geq 0}$ and satisfies 
\begin{equation}
\label{eq-icphi}
\mu\PAR{\frac{\Gamma_2(\Phi'(g))}{\Phi''(g)}}\geq  \rho \, \mu\PAR{\frac{\Gamma(\Phi'(g))}{\Phi''(g)}}
\end{equation}
for all $I$-valued functions $g$, then $\mu$ satisfies  the $\Phi$-entropy inequality~\eqref{PHII} for all $I$-valued functions~$f$. 
\end{eprop}


Let us first note that any admissible function $\Phi$ is stricly convex, so that its derivative $\Phi'$ is increasing and has an inverse $\Phi'^{-1}:$ then the integral criterion~\eqref{eq-icphi} writes
$$
\mu\PAR{\frac{\Gamma_2(h)}{\Phi'' \circ \Phi'^{-1} (h)}}\geq  \rho \, \mu\PAR{\frac{\Gamma(h)}{\Phi'' \circ \Phi'^{-1} (h)}}
$$
for all $I$-valued functions $h$ with values in the image of $\Phi'.$
 This integral criterion appears in \cite[proof of Theorem 7.2.2]{helffer02} in the case when $L = \Delta - < \nabla V,  \nabla >$ and $\mu$ is the reversible ergodic measure $e^{-V}$. It extends the criteria~\eqref{critereint} for the Poincar\'e inequality, with $\Phi(x) = x^2$ and $\Phi'' \circ \Phi'^{-1} (x) =2$, and~\eqref{criteresuperint} for the logarithmic Sobolev inequality, with $\Phi(x) = x \ln x$ and $\Phi'' \circ \Phi'^{-1} (x) =e^{1-x}.$ Let us point out that in this diffusion setting {\it the measure $\mu$ need not be assumed to be reversible}, as it is usually the case in the previous works (see \cite[Proposition 5.5.6]{logsob} for instance).
 
 For the $\Phi_p$ maps with $p$ in $]1,2[$ it writes
\begin{equation}
\label{eq-icphipp}
\mu\PAR{g^{\frac{2-p}{p-1}}\Gamma_2(g)}\geq  \rho \, \mu\PAR{g^{\frac{2-p}{p-1}}\Gamma(g)}
\end{equation}
for all positive functions $g.$

\bigskip

Let us note that this family of Beckner's inequalities for $p$ in $]1,2[$ has been obtained in \cite{dns} under a different integral condition, which does not seem to be comparable to our condition~\eqref{eq-icphipp} for $p$ in $]1,2[.$

\begin{erem}\label{remcritereint}
At least for $p$ close to $1$ the integral criterion~\eqref{eq-icphipp} is not equivalent to~\eqref{beckneravecp}-\eqref{beckner}  (hence strictly stronger), thus extending the case of the logarithmic Sobolev inequality when $p=1:$

Following B. Helffer (see \cite[p.~114]{helffer02} or \cite[p. 91]{logsob}) we build a probability measure $\mu$ on $\rr$ such that~\eqref{beckner} holds while~\eqref{eq-icphipp} does not hold. We consider the generator $Lf=f''-\Psi'f'$ on $\rr$, with $\Psi (x) =x^4-b x^2$, and its reversible ergodic measure $\mu(dx)=\exp(-\Psi(x))dx/Z$ where $Z$ is a normalization constant. For any $b$ the measure $\mu$ satisfies a logarithmic Sobolev inequality, hence the Beckner inequality~\eqref{beckner}. But, letting $g(x)=\exp(-x^2)$ and $b=1+p/(p-1)$, we obtain 
$$
\mu\PAR{g^{\frac{2-p}{p-1}}\Gamma_2(g)}= \int \PAR{\PAR{4x^2 - 2}^2+48x^4}\frac{e^{x^2-x^4}}{Z}dx-8 \, \PAR{1+\frac{p}{p-1}}\int x^2\frac{e^{x^2-x^4}}{Z}dx
$$   
which is negative for  $p$ in $]1,2[$ close to 1. Hence~\eqref{eq-icphipp} cannot hold since the right hand side in nonnegative, so that, for these $p$, the integral criterion~\eqref{eq-icphipp} is not a necessary condition to the $\Phi_p$-entropy inequality~\eqref{beckneravecp}-\eqref{beckner}.
\end{erem}

\noindent\emph{\textbf{Proof of Proposition~\ref{prop-ci}}}

\proofbegin\,
The argument follows the argument of Theorem~\ref{theophisobolev}.  If $f$ is a given positive function we let
$$
H(u)=\mu\PAR{\Phi(\PT{u} f)}, \quad u \geq 0.
$$

For $t >u$ fixed we let again $\psi (s) =  \PT{s} (\Phi(\PT{t-s} f))$, so that
$$
H(u) = \mu \big(   \PT{t-u} (\Phi(\PT{u} f)) \big) = \mu ( \psi(t-u) )
$$
since the ergodic measure $\mu$ is necessarily invariant.  In particular  Lemma~\ref{lemmepsi'} ensures that
 \begin{equation}\label{H'}
H'(u)=- \mu(\psi'(t-u)) = - \mu \big(  \Phi''(\PT{u} f) \Gamma(\PT{u} f)  \big)= - \mu\PAR{\frac{\Gamma(\Phi'(\PT{u} f))}{\Phi''(\PT{u} f)}} 
\end{equation}
and 
\begin{equation}\label{H''}
H''(u)=2 \mu\PAR{\frac{\Gamma_2(\Phi'(\PT{u} f))}{\Phi''(\Pt{u} f)}}
+
\mu \PAR{\PAR{\frac{\Gamma (\Phi'(\PT{u}   f))}{\Phi''(\PT{u}   f)}}^2\PAR{\frac{-1}{\Phi''}}''(\PT{u}f)}
\end{equation}
again by the invariance property of $\mu.$ 

The second term on the right hand side of~\eqref{H''} is nonnegative by assumption on $\Phi$, so that
$$
H''(u) \geq - 2\, \rho \, H'(u), \quad u \geq 0
$$
by the integral criterion~\eqref{eq-icphi}. Integrating between $0$ and $t$ gives
$$
-H'(t) \leq -H'(0) \, e^{-2 \rho t}
$$
and integrating again between $0$ and $+\infty$ concludes the argument by ergodicity of $(\Pt)_{t \geq 0}$.
\proofend


\subsection{Large time behaviour of the Markov semigroup and the linear Fokker-Planck equation}\label{subsectFP}

We now turn to the long time behaviour of the diffusion Markov semigroup and of the solutions to an associated linear Fokker-Planck equation.

\medskip

The argument for the Markov semigroup is simpler:  Let  $(\Pt )_{t\geq 0}$ be a diffusion semigroup, ergodic for the probability measure $\mu$.  If $\Phi$ is a $\mathcal C^2$ strictly  convex function on an interval $I$, then \eqref{H'} above writes
\begin{equation}\label{HI}
\frac{d}{dt}\entf{\mu}^\Phi\PAR{\Pt f}=-\mu\PAR{\frac{\Gamma(\Phi'(\Pt f))}{\Phi''(\Pt f)}}
\end{equation}
for all $t \geq 0$ and all $I$-valued functions $f.$ If $C$ is a positive number, then there is equivalence between:
\begin{enumerate}[(i)]
\item the measure $\mu$ satisfies the $\Phi$-entropy inequality 
\begin{equation}
\label{eq-phisob}
\entf{\mu}^\Phi\PAR{f}\leq C\mu \PAR{\frac{\Gamma(\Phi'( f))}{\Phi''( f)}}
\end{equation}
 for all $I$-valued functions $f$;
\item the semigroup converges in $\Phi$-entropy with exponential rate: 
\begin{equation}
\label{cv1}
\entf{\mu}^\Phi\PAR{\Pt f}\leq e^{-\frac{t}{C}}\entf{\mu}^\Phi\PAR{ f}
\end{equation}
for all $t \geq 0$ and all $I$-valued functions $f.$
\end{enumerate}

Indeed (ii) follows from (i) by \eqref{HI} and (i) follows from (ii) by differentiation at $t=0.$

\bigskip

We now turn to the study of the linear Fokker-Planck equation
\begin{equation}
\label{eq-fp1}
\frac{\partial u_t}{\partial t}  = \textrm{div} \SBRA{  D(x)( \nabla u_t + u_t (\nabla V(x)+F(x)))}, \quad t \geq 0, \, x \in \rr^n
\end{equation}
where div stands for the divergence, $D(x)$ is a positive symmetric matrix $\dR^n$ and the vector field $F$ satisfies the condition
\begin{equation}\label{rev}
\textrm{div} \PAR{e^{-V} D F}=0.
\end{equation}

It is one of the purposes of~\cite{amtucpde01} and~\cite{arnoldcarlenju08} to rigorously study the asymptotic behaviour of solutions to~\eqref{eq-fp1}-\eqref{rev}.  Let us formally rephrase the argument in our semigroup terminology.

 We let $L$ be the Markov diffusion generator defined by 
 \begin{equation}\label{fparnold2}
Lf  =  \textrm{div} (D \nabla  f)  -  < D ( \nabla V-F) , \nabla f>.
\end{equation}
Let us assume that $L$ satisfies the $CD(\rho,\infty)$ criterion with $\rho>0$, that is, for instance in the case when $D(x) = d(x) I$ is a scalar matrix and $F=0$,
\begin{equation}
\label{cdrhoarnold}
\big(\frac{1}{2} - \frac{n}{4} \big) \frac{1}{d} \nabla d \otimes \nabla d + \frac{1}{2} ( \Delta d \, - <\nabla d , \nabla V > ) I + d \, {\textrm Hess} V 
+ \frac{1}{2} (\nabla V \otimes \nabla d + \nabla d \otimes \nabla V) - \textrm{Hess} \, d \geq \rho \, I.
\end{equation}
Then the semigroup $(\Pt)_{t \geq 0}$ associated to $L$ is ergodic and the ergodic probability measure is explicitely given by $d\mu=e^{-V}/Zdx$ where $Z$ is a normalization constant. If moreover $\Phi$ is an admissible function on an interval $I$, then the $\Phi$-entropy inequality~\eqref{eq-phisob} holds with  $C=1/(2\rho)$ by Theorem~\ref{theophisobolev}, so that the semigroup converges to $\mu$ according to~\eqref{cv1}. 

But, under  the condition~\eqref{rev},  the solution to the Fokker-Planck equation~\eqref{eq-fp1} for the initial datum $u_0$ is given by $u_t=e^{-V}\,\Pt (e^V u_0)$. Then we can deduce the convergence of the solution $u_t$ towards the stationary state $e^{-V}$ (up to a multiplicative constant) from the convergence estimate~\eqref{cv1} for the Markov semigroup, in the form
$$
\entf{\mu}^\Phi\PAR{\frac{u_t}{e^{-V}}}\leq e^{-2\rho{t}}\entf{\mu}^\Phi\PAR{ \frac{u_0}{e^{-V}}}, \quad t \geq 0
$$
for all initial data $u_0$ such that the map $e^V u_0$ be $I$-valued. 

\bigskip

In fact one can obtain estimates on the long time behaviour of solutions to~\eqref{eq-fp1}  {\it without the condition~\eqref{rev}}. Let us indeed consider the linear Fokker-Planck equation
\begin{equation}
\label{eq-fp2}
\frac{\partial u_t}{\partial t}  = \textrm{div} \SBRA{  D(x) (\nabla u_t + u_t a(x) )}, \quad t \geq 0, \, x \in \rr^n
\end{equation}
where again $D(x)$ is a positive symmetric matrix $\dR^n$ and $a(x) \in \rr^n.$ Its generator is the dual $L^*$ (with respect to the Lebesgue measure) of the generator 
\begin{equation}\label{exam}
Lf = \textrm{div} (D \nabla  f)  -  < D a , \nabla f>. 
\end{equation}
Assume that the semigroup associated to $L$ is ergodic and that its invariant probability measure $\mu$ satisfies a $\Phi$-entropy inequality~\eqref{eq-phisob} with a constant $C\geq0$: this holds for instance if $L$ satisfies the $CD(1/(2C),\infty)$ criterion, that is, if
$$
D Ja(x) D +(D Ja(x) D)^*\geq \frac{1}{C} D
$$
for all $x\in\dR^n$ if $D$ is constant (by~\eqref{cdrhocst}), or if \eqref{cdrhoscalaire} holds if $D(x)$ is a scalar matrix, and so on.

In this more general setting when $a (x)$ is not the gradient of a potential, the invariant measure is not explicit. Moreover the explicit relation between the solution $u_t$ of the linear Fokker-Planck and the associated semigroup $\Pt $ does not hold anymore,  which lead above from the asymptotic behaviour of the semigroup to that of the solutions of the Fokker-Planck equation. 

This can be replaced by the following argument, for which we only assume that the ergodic measure has a positive density $u_{\infty}$ with respect to the Lebesgue measure:

 Let $u$ be a solution of~\eqref{eq-fp2} for the initial datum $u_0.$ Then
$$
\frac{d}{dt}\entf{\mu}^\Phi\PAR{\frac{u_t}{u_\infty}}=\int \Phi'  \PAR{\frac{u_t}{u_\infty}}L^* u_t \, dx=\int  L \Big[  \Phi'  \PAR{\frac{u_t}{u_\infty}} \Big] \, \frac{u_t}{u_\infty} d\mu
=
-\int {\Phi''\PAR{\frac{u_t}{u_\infty}}}\Gamma\PAR{{\frac{u_t}{u_\infty}}}\,d\mu
$$
by Lemma~\ref{end-end} with  $f=\frac{u_t}{u_\infty}$ and $\phi = \Phi'.$ Then the $\Phi$-Entropy inequality~\eqref{eq-phisob} for $\mu$ implies the exponential convergence
$$
\entf{\mu}^\Phi\PAR{\frac{u_t}{u_\infty}}\leq e^{-\frac{t}{C}}\entf{\mu}^\Phi\PAR{ \frac{u_0}{u_\infty}}, \quad t \geq 0.
$$

\begin{elem}
\label{end-end}
 Let $L$ be a diffusion generator with invariant measure~$\mu$.  Then
$$
\int  L{ \phi  \PAR{f}}\, fd\mu =\int  Lf\,{ \phi  \PAR{f}}d\mu= -\int\Gamma\PAR{{f,\phi  \PAR{f}}}\,d\mu
$$
for all functions $f$ and all one-to-one functions $\phi$.
\end{elem}
\begin{eproof}
Let $g=\phi(f)$, so that
$$
\int  L{ \phi \PAR{f}}\, fd\mu =\int Lg\, \psi(g)d\mu,
$$
where $\psi=\phi^{-1}$. Then if $\Psi$ is an antiderivative of $\psi$ we get
$$
\int Lg\, \psi(g)d\mu=\int\PAR{ Lg\, \psi(g)-L\Psi(g)}d\mu = -\int \Psi''(g)\Gamma(g)d\mu=-\int \Gamma(f,g)d\mu,
$$
by invariance of $\mu$ and the diffusion properties~\eqref{eq-diffusion}. This concludes the argument by the identity $\psi'(\phi(x))\phi'(x)=1$.
\end{eproof}

\bigskip

Hence (for instance) we have formally obtained:

\begin{theo}\label{theo-cvfp}
In the above notation, let $\Phi$ an admissible function and assume that the generator $L$ of~\eqref{exam} satisfies the integral criterion~\eqref{eq-icphi} to a $\Phi$-entropy inequality and has an ergodic measure with smooth positive density $u_{\infty}.$ Then all solutions $u = (u_t)_{t \geq 0}$ to the Fokker-Planck equation \eqref{eq-fp2} converge to $u_{\infty}$ in $\Phi$-entropy, with
$$
\entf{\mu}^\Phi\PAR{\frac{u_t}{u_\infty}}\leq e^{-2 \rho t}\entf{\mu}^\Phi\PAR{ \frac{u_0}{u_\infty}}, \quad t \geq 0.
$$
\end{theo}

\bigskip

The three sections below are devoted to improvements and extensions of the $\Phi$-entropy inequalities considered in this section. First of all, in Sections \ref{sect12} and \ref{sectnot12}, we improve Theorem \ref{theophisobolev} for $\Phi_p$-entropies which are the main examples of such $\Phi$-entropies. In Section~\ref{sectnot12} we also derive upper and lower bounds on $\entf{\Pt}^\Phi (f)$, analogous to \eqref{localPHII} and \eqref{reverselocalPHII} and still equivalent to the $CD(\rho, \infty)$ criterion, for other maps $\Phi$ that do not satisfy the admissibility hypotheses of Theorem \ref{theophisobolev}.  In a last section we shall see how these results extend to the setting of inhomogeneous Markov semigroups and to a instance of nonlinear evolution problem.

\section{Refined $\Phi$-entropy inequalities in the admissible case}
\label{sect12}


In this section we give and study improvements of Theorem~\ref{theophisobolev} for the main family of admissible $\Phi$, namely the $\Phi_p$ functions given by~\eqref{eq-phip}, for $p \in ]1,2[.$ 

\begin{theo}
\label{thm-main}
Let $\rho$ be a real number and  $p$ in $]1,2[$. Then the following assertions are equivalent, with $\disp  \frac{1 - e^{-2\rho t}}{\rho}$ and $\disp \frac{e^{2\rho t} -1}{\rho}$ replaced by $2 t$ if $\rho =0$:
\begin{enumerate}[(i)]
\item
 the semigroup $(\PT{t})_{t \geq 0}$ satisfies the ${CD}(\rho,\infty)$ criterion;
\item 
 the semigroup $(\PT{t})_{t \geq 0}$ satisfies the refined local $\Phi_p$-entropy inequality
 \begin{equation}\label{refinedlocalPHII}
\frac{1}{(p-1)^2}\SBRA{{\Pt(f^p)}-{\Pt(f)^p}\PAR{\frac{\Pt(f^p)}{\Pt(f)^p}}^{\frac{2}{p}-1}}\leq \frac{1-e^{-2\rho t}}{\rho}\Pt\PAR{f^{p-2} \, {\Gamma(f)}}
\end{equation}
 for all positive $t$ and all positive functions $f$;
 \item
the semigroup $(\PT{t})_{t \geq 0}$ satisfies the reverse local $\Phi$-entropy inequality  
$$
\frac{1}{(p-1)^2}\SBRA{{\Pt(f^p)}-{\Pt(f)^p}\PAR{\frac{\Pt(f^p)}{\Pt(f)^p}}^{\frac{2}{p}-1}}\geq
 \disp\frac{e^{2\rho t}-1}{\rho}\PAR{\frac{(\Pt f)^p}{\Pt\PAR{f^p}}}^{\frac{2}{p}-1}\PAR{\Pt f}^{p-2} \, {\Gamma(\Pt f)}
$$
 for all positive $t$ and all positive functions $f$.
  \end{enumerate}
\end{theo}

\begin{eproof}
We first assume that $(i)$ holds and prove $(ii)$ and $(iii).$ As in the proof of Theorem~\ref{theophisobolev} we let $\psi (s) = \PT{s}\PAR{\Phi_p\PAR{\PT{t-s} f}}$, where in this proof and in the proofs  of Proposition~\ref{prop-cip} only $\Phi_p (x) = x^p/(p(p-1));$ then Lemma \ref{lemmepsi'} specifies as
$$
\psi''(s)=2 \, \PT{s}\PAR{\frac{\Gamma_2(\Phi_p'(\PT{t-s}   f))}{\Phi''_p(\PT{t-s}   f)}}
+ (2-p)(p-1)\, \PT{s}\PAR{\PAR{\frac{\Gamma (\Phi'_p(\PT{t-s}   f))}{\Phi''_p(\PT{t-s}   f)}}^2 \frac{1}{(\PT{t-s}f)^p}}.
$$

The second term on the right hand side is nonnegative for $p$ in $]1,2[$. By the Cauchy-Schwarz inequality it is even bounded by below by
$$
\PAR{\PT{s}\PAR{\frac{\Gamma (\Phi'_p(\PT{t-s}   f))}{\Phi''_p(\PT{t-s}   f)}}}^2 \frac{1}{\PT{s}\PAR{(\PT{t-s}f)^p}} = \frac{1}{p(p-1)} \frac{\psi'(s)^2}{\psi(s)}.
$$
Then the $CD(\rho,\infty)$ criterion implies 
$$
\psi''(s)\geq 2\rho\psi'(s)+\frac{2-p}{p}\frac{\psi'(s)^2}{\psi(s)}, \quad s \in [0,t],
$$
that is, 
$$
\PAR{\frac{\psi'(s)}{\psi(s)^{(2-p)/p}}e^{-2\rho s}}' \geq 0, \quad  s\in[0,t].
$$ 
For $0\leq u\leq v\leq t$ we integrate over the interval $[u,v]$ to obtain 
$$
\frac{\psi'(v)}{\psi(v)^{(2-p)/p}}e^{-2\rho v}\geq \frac{\psi'(u)}{\psi(u)^{(2-p)/p}}e^{-2\rho u}.
$$

For $u=s$ and $v=t$, integrating in $s$ over the set $[0,t]$ yields 
$$
\frac{p}{p-1} \big[ \psi(t)^{2(p-1)/p} - \psi(0)^{2(p-1)/p} \big] \leq\frac{\psi'(t)}{\psi(t)^{(2-p)/p}}\frac{1-e^{-2\rho t}}{\rho},
$$
which leads to $(ii)$.

 For $u=0$ and $v=s$, integrating over the set $[0,t]$ yields 
$$
\frac{p}{p-1} \big[ \psi(t)^{2(p-1)/p} -  \psi(0)^{2(p-1)/p} \big] \geq \frac{\psi'(0)}{\psi(0)^{(2-p)/p}}\frac{e^{2\rho t}-1}{\rho},
$$
which leads to $(iii)$. 

\medskip

Let us conversely assume that $(ii)$ or $(iii)$ holds. For $f=1+\varepsilon g$, the left hand side in $(ii)$ and $(iii)$  is 
$$
\varepsilon^2 \, \var{\Pt}{g}+o(\varepsilon^2)
$$
and the right hand side is
$$
\varepsilon^2\frac{1-e^{-2\rho t}}{\rho}  \, \Pt\PAR{\Gamma(g)} + o(\varepsilon^2) \quad \text{and} \quad 
\varepsilon^2\frac{e^{2\rho t} - 1}{\rho}  \, \Pt\PAR{\Gamma(g)} + o(\varepsilon^2)
$$
respectively. Hence, as $\varepsilon$ goes to 0, $(ii)$ implies the local Poincar\'e inequality~\eqref{localPI} whereas $(iii)$ implies the reverse local Poincar\'e inequality~\eqref{reverselocalPI}, which are both equivalent to the $CD(\rho,\infty)$ criterion. This concludes the proof of the theorem.    
\end{eproof}

The heat equation on $\dR^n$ satisfies the $CD(0,\infty)$ criterion and is linked with the standard Gaussian measure $\gamma$ by the identity $\PT{1/2} \,  g(0) = \gamma (g).$ Hence, applying ~\eqref{refinedlocalPHII} with $\rho =0$  to this semigroup at $t=1/2$ and $x=0$ leads to the following bound for the Gaussian measure:
$$
\frac{1}{(p-1)^2}\SBRA{{\gamma(f^p)}-{\gamma(f)^p}\PAR{\frac{\gamma(f^p)}{\gamma(f)^p}}^{\frac{2}{p}-1}}\leq \gamma\PAR{f^{p-2} \, \Gamma(f)}.
$$

The Gaussian measure is also the ergodic measure of the Ornstein-Uhlenbeck semigroup, and for ergodic measures we obtain the more general result:
 
\begin{ecor}\label{coro-ad}
In the above notation, if the semigroup $(\Pt)_{t \geq 0}$ is $\mu$-ergodic and sat\nobreak isfies the $CD(\rho, \infty)$ criterion with $\rho>0$, then the measure $\mu$ satisfies 
\begin{equation}
\label{eq-ad}
\frac{1}{(p-1)^2}\SBRA{{\mu(f^p)}-{\mu(f)^p}\PAR{\frac{\mu(f^p)}{\mu(f)^p}}^{\frac{2}{p}-1}}\leq \frac{1}{\rho}\mu\PAR{f^{p-2} \, \Gamma(f)}
\end{equation}
for all positive functions $f$ or equivalently 
\begin{equation}
\label{eq-ad2}
\frac{p^2}{(p-1)^2}\SBRA{{\mu(g^2)}-{\mu(g^{2/p})^p}\PAR{\frac{\mu(g^2)}{\mu(g^{2/p})^p}}^{\frac{2}{p}-1}}\leq \frac{4}{\rho}\mu\PAR{\Gamma(g)}
\end{equation}
for all positive functions $g$. 
\end{ecor}

The refined $\Phi_p$-entropy inequality~\eqref{eq-ad} has been obtained by A. Arnold and J.~Dolbeault in \cite{arnolddolbeault05}  for the generator $L$ defined in \eqref{fparnold2} with $D(x)$ a scalar matrix and $F=0$  and for the ergodic measure $\mu = e^{-V}$, and under the corresponding $CD(\rho, \infty)$ condition \eqref{cdrhoarnold}.

\bigskip

As pointed out in \cite{arnolddolbeault05}, the bound given by Corollary \ref{coro-ad} improves on the Beckner inequality
$$
\frac{\mu(g^2) - \mu(g^{2/p})^p}{p-1} \leq \frac{2}{p   \rho} \, \mu (\Gamma (g))
$$
given by Theorem \ref{thm-main} since 
\begin{equation}\label{comp-beckner+}
\frac{\mu(g^2) - \mu(g^{2/p})^p}{p-1} \leq \frac{p}{2 (p-1)^2} \Big[ \mu(g^2) -\mu(g^{2/p})^p \, \Big(\frac{\mu(g^2)}{\mu(g^{2/p})^p}\Big)^{\frac{2}{p}-1}  \Big]
\end{equation}
for all positive functions $g$.

\bigskip
In the first section we have noticed that for all $g$ the map 
$$
p \mapsto \left\{
\begin{array}{cl}
\disp \frac{\mu(g^2) - \mu(g^{2/p})^p}{p-1} & \text{ if }p \neq 1\\
\disp \ent{\mu}{g^2} & \text{ if }p=1
\end{array}
\right.
$$
is continuous and nonincreasing on $]0, +\infty[$
and takes the value $\var{\mu}{g}$ at $p=2$.

In the next proposition we show similar properties for the functional introduced in \eqref{comp-beckner+}:

\begin{prop}
\label{prop-7}
For any Borel probability measure $\mu$ on $\rr^n$ and any positive  $g$ on $\rr^n$ the map 
$$
p \mapsto \frac{p}{2 (p-1)^2} \Big[ \mu(g^2) -\mu(g^{2/p})^p \, \Big(\frac{\mu(g^2)}{\mu(g^{2/p})^p}\Big)^{\frac{2}{p}-1}  \Big]
$$ 
is nonincreasing on $]1,+\infty[.$ Moreover it takes the value $\var{\mu}{g}$ at $p=2$ and admits the limit $ \ent{\mu}{g^2}$ as $p$ tends to $1.$

\end{prop}

\begin{erem}
The map is also continuous on the left hand side of $1$, but not necessarily monotone on $]0,1[$: for instance, if $\mu$ is the standard 
Gaussian measure on $\rr$ and $f$ the map defined in $\rr$ by $f(x) = \disp \sqrt{2} e^{-x^2/2}$, then the map takes the approximate values $0,061$ at $p=0,1$, then $0,134$ at $p=0,5$ and $0,103$ at $p=0,9.$
\end{erem}

\begin{erem}
\label{lsibeckner+}
If a measure $\mu$ satisfies the logarithmic Sobolev inequality with constant $C$
$$
\ent{\mu}{g^2} \leq C \, \mu ( \Gamma(g))
$$
for all functions $g$, then it follows from Proposition \ref{prop-7} that for all $p$ in $]1,2[$ it satisfies the refined $\Phi_p$-entropy inequality \eqref{eq-ad2} with the constant $2 p C$ instead of $4/\rho.$ 

For instance, if the  $CD(\rho,\infty)$ criterion with $\rho>0$ holds for an $\mu$-ergodic semigroup $(\Pt)_{t\geq 0}$ then we can take $C = 2/\rho$ and we recover that for all $p>1$ the measure $\mu$ satisfies~\eqref{eq-ad2}, but  with a constant $4 p/\rho$ instead of the finer constant $4/\rho$ given by Corollary~\ref{coro-ad}.

In the same way, assume that the measure $\mu$ satisfies the refined $\Phi_p$-entropy inequality \eqref{eq-ad2} for a $p$ $]1, 2]$ and a constant $4/\rho,$ and let $q$ in $[p,2].$ Then, by Proposition~\ref{prop-7},  the measure $\mu$ satisfies~\eqref{eq-ad2} with $q$ instead of $p$ and a constant $4 q / p \rho$, instead of the finer constant $4/\rho$ given by Corollary~\ref{coro-ad} under the $CD(\rho,\infty)$ criterion.
\end{erem}

\noindent\emph{\textbf{Proof of Proposition~\ref{prop-7}}}

\proofbegin\,
We first prove that the map is nonincreasing in $p.$ By homogeneity we may assume that $\mu(g^2) =1$, and in the notation $h = \ln g^2$ and $t = 1/p$ we prove that for any function $h$ such that $\mu(e^h) =1$ the map
$$
t \mapsto \disp \frac{t}{(t-1)^2} \Big[ 1 - \mu(e^{th})^{2(1/t-1)}\Big]
$$
is nondecreasing on $]0, 1[.$

For this purpose we let fixed such an $h$ and we prove that the map $a$ defined on $]0,1[$ by
$$
a(t)=\frac{t}{(1-t)^2}(1-e^{b(t)})
$$ 
where $b(t)=\frac{2}{t}(1-t)\ln \mu(e^{th}).$

Its derivative is
$$
a'(t) = \disp \frac{t+1}{(1-t)^3} e^{b(t)} \Big[ e^{-b(t)} - \big( 1 + \frac{t(1-t)}{1+t} \, b'(t) \big) \Big]
$$
so it is sufficient to prove that
 $$
 e^{-b(t)} \geq 1 + \frac{t(1-t)}{1+t} \, b'(t).
 $$
But $e^x \geq 1 + x + x^2/2$ for all $x \geq 0$ and $-b(t) \geq 0$ since $\ln \mu(e^{th}) \leq \ln (\mu(e^h)^t) = 0$ by the H\"older inequality, so it is sufficient to prove that
$$
-b (t)+\frac{b(t)^2}{2}\geq \frac{t(1-t)}{1+t} \, b'(t)
$$
or equivalently that
$$
c(t)=-\ln \mu( e^{th}) +\frac{1-t^2}{t^2}\mu( e^{th})^2-(1-t)\frac{\mu( he^{th})}{\mu( e^{th})}
$$
is nonnegative. 

For those $t$ (if any) such that $d(t) = \mu (h \, e^{th})$ is nonpositive, then $c(t)$ is nonnegative as the sum of nonnegative terms.

Morever $d$ is nondecreasing since $d'(t)= \mu(h^2 \, e^{th}) \geq 0$. Hence, if there exists $t_0 \in ]0,1[$ such that $d(t_0)>0, $ then $d(t) >0$ on an interval $]t_1, 1[:$ on $]0, t_1]$ we have $d(t) \leq 0$ so that $c(t) \geq 0$ and remains to be proven that $c(t) \geq 0$ on $]t_1, 1[$ on which $d(t) >0.$

Now
$$
c'(t) = - \frac{2}{t^3} \ln^2 \mu(e^{th}) + 2 \frac{1-t^2}{t^2} \, \ln \mu(e^{th}) \, \frac{\mu(h \, e^{th})}{\mu(e^{th})}\\
 +(1-t) \frac{ \mu(h \, e^{th})^2-\mu(h^2 \, e^{th}) \, \mu(e^{th}) }{\mu(e^{th})^2}
$$
where the first term on the right hand side is nonpositive for all $t$, the second term nonpositive for all $t$ in $]t_1, 1[$ and the third term nonpositive for all $t$ by the Cauchy-Schwarz inequality. Hence $c$ is a nonincreasing function on $]t_1, 1[,$ and $c(t) \geq \lim_{s \to 1} c(s)  =0$ for all $t$ in $]t_1, 1[.$

As a consequence $c$ is nonegative on $]0, 1[$ and $a$ is indeed nonincreasing on $]0, 1[.$

\bigskip
Then we prove that the functional admits the limit $\ent{\mu}{g^2}$ as $p$ tends to $1$. By homogeneity we may again assume that $\mu(g^2) =1$, and we let $t = 1/p$ and $h = \ln g^2,$ so that $\mu(e^h) =1$. Then we have to prove that
$$
\frac{t}{2 (t-1)^2} \Big[ 1 - \mu(e^{th})^{2(1/t-1)}\Big]
$$
tends to $\mu(h \, e^h)$ as $t$ tends to $1$, which can be checked by letting $t=1+ \varepsilon$ and performing Taylor expansions around $\varepsilon = 0$.
\proofend

\bigskip

We conclude this section by proving that the integral criterion~\eqref {eq-icphipp} for the $\Phi_p$-entropy inequality~\eqref{beckneravecp}-\eqref{beckner} of Section~\ref{sub-11} is also a sufficient condition for  the stronger inequality~\eqref{eq-ad}-\eqref{eq-ad2}:

\begin{eprop}
\label{prop-cip}
Let $\rho$ be a positive number and $p$ in $]1,2[$. If the probability measure $\mu$ is ergodic for the diffusion semigroup $(\Pt)_{t \geq 0}$ and satisfies 
$$
\mu\PAR{g^{\frac{2-p}{p-1}}\Gamma_2(g)}\geq  \rho \, \mu\PAR{g^{\frac{2-p}{p-1}}\Gamma(g)}
$$
for all positive functions $g$, then $\mu$ satisfies  the refined $\Phi_p$-entropy inequality~\eqref{eq-ad} for all positive functions $f$. 
\end{eprop}

\begin{erem} 
As pointed out in Remark~\ref{remcritereint},  the integral criterion~\eqref{eq-icphipp} of Proposition~\ref{prop-cip} is not equivalent (hence strictly stronger) to the Beckner inequality~\eqref{beckneravecp}-\eqref{beckner} at least for $p$ close to $1.$ It is not either equivalent to the stronger inequality~\eqref{eq-ad}-\eqref{eq-ad2}: for the counter-example of Remark~\ref{remcritereint}, the criterion does not hold, although the considered measure $\mu$ satisfies a logarithmic Sobolev inequality, hence~\eqref{eq-ad}-\eqref{eq-ad2} by Remark~\ref{lsibeckner+}. 
\end{erem}

\noindent\emph{\textbf{Proof of Proposition~\ref{prop-cip}}}

\proofbegin\,
The argument follows the argument of Proposition~\ref{prop-ci} by taking advantage of elements of the proof of Theorem~\ref{thm-main}.  If $f$ is a given positive function, the function
$$
H(u)=\mu\PAR{\Phi_p(\PT{u} f)}, \quad u \geq 0,
$$
where again $\Phi(p) = x^p/(p(p-1))$, has a second derivative given by
$$
H''(u) = 2 \, \mu\PAR{\frac{\Gamma_2(\Phi_p'(\PT{u} f))}{\Phi_p''(\PT{u} f)}}
+ (2-p)(p-1) \, \mu\PAR{\frac{1}{(\PT{u} f)^p}\PAR{\frac{\Gamma(\Phi_p'(\PT{u} f))}{\Phi_p''(\PT{u} f)}}^2}.
$$
Hence it satisfies the inequality
\begin{equation}\label{H'H''}
H''(u) \geq -2 \, \rho \, H'(u)-\frac{p-2}{p}\frac{H'(u)^2}{H(u)}, \quad u \geq 0
\end{equation}
by the integral criterion and the Cauchy-Schwarz inequality. For all $t\geq 0$, integrating over the set $[0,t]$ gives
$$
H'(t) H(t)^{(p-2)/p} \geq H'(0) H(0)^{(p-2)/p} \, e^{-2\rho t}.
$$
Integrating between $0$ and $+\infty$ conclude the argument by ergodicity of $(\Pt)_{t \geq 0}$. 
\proofend

 \medskip
 
 \begin{erem}
 For $\rho=0$, and following \cite{arnolddolbeault05}, the convergence of $\Pt f$ towards $\mu(f)$ can be measured as
$$
 \vert H'(t) \vert \leq \frac{\vert H'(0)\vert}{1+ \alpha t}, \quad t \geq 0
 $$
 where $\alpha = \frac{2-p}{p} \vert H'(0) \vert / \entf{\mu}^{\Phi_p}(f).$ This is an illustration of the improvement given by using the second term on the right hand side in \eqref{H'H''}. However it does not give a rate of convergence to $0$ on $\entf{\mu}^{\Phi_p}(\Pt f)$, as in Section~\ref{subsectFP} for $\rho >0.$ Indeed, for $\rho =0,$ equation \eqref{H'H''} is equivalent to $H^{2-2/p}$ being convex: hence \eqref{H'H''} is solved by positive maps decaying to $0$ at infinity as slowly as we like.

 \end{erem}

\section{The case of a non-admissible function}
\label{sectnot12}

We now turn to the issue of deriving $\Phi$-entropy inequalities for maps $\Phi$ which do not satisfy the admissibility assumptions of Theorem~\ref{theophisobolev}.

As shown in Section~\ref{sectphi}, for a given probability measure $\mu$ and a positive function $g$, the map
$$
p \mapsto \left\{
\begin{array}{cl}
\disp \frac{\mu(g^2) - \mu(g^{2/p})^p}{p-1} & \text{ if }p \neq 1\\
\disp \ent{\mu}{g^2} & \text{ if }p=1
\end{array}
\right.
$$
is  nonincreasing on $]0, +\infty[$. Hence, if a measure $\mu$ satisfies the logarithmic Sobolev inequality with constant $C$
$$
\ent{\mu}{g^2} \leq C \, \mu ( \Gamma(g))
$$
for all functions $g$, then for all $p >1$ it satisfies the Beckner inequality
$$
\frac{\mu(g^2) - \mu(g^{2/p})^p}{p-1}\leq C \mu ( \Gamma(g))  
$$
for all positive functions $g$ with the same constant $C$, 
or equivalently 
\begin{equation}
\label{eq-phipp}
\frac{\mu(f^p) - \mu(f)^p}{p-1}\leq C\frac{p^2}{4} \mu( f^{p-2}\Gamma(f))  
\end{equation}
for all positive functions $f.$ 

For instance, if the $CD(\rho,\infty)$ criterion with $\rho>0$ holds for an $\mu$-ergodic semigroup $(\Pt)_{t\geq 0}$ then we can take $C = 2/\rho$ and for all $p>1$ the measure $\mu$ satisfies~\eqref{eq-phipp} with $C p^2 /4 = p^2 /(2 \rho).$

\medskip

For $1 < p \leq 2$, and under the $CD(\rho, \infty)$ criterion, Theorem \ref{theophisobolev} ensures that inequality~\eqref{eq-phipp} even holds for the semigroup $(\PT{t})_{t\geq 0}$ in the local form

\begin{equation}\label{becknerloc}
\frac{{\PT{s}(f^p)}-\PT{s}(f)^p}{p-1} \leq p \frac{1 - e^{-2 \rho s}}{2 \rho} \, \PT{s}(f^{p-2} \Gamma (f)), \quad s \geq0
\end{equation}
(note the finer constant $p / (2 \rho)$ instead of $p^2 / (2 \rho)$ in the right hand side).

\medskip

For $p >2$ this may not hold.
Indeed, for $f=\PT{t-s} g$ with $t$ fixed, \eqref{becknerloc} writes 
$$
\psi(s)-\psi(0)\leq \psi'(s)\frac{1-e^{-2\rho s}}{2\rho}
$$
in the notation $\psi\PAR{s} =  \PT{s}\PAR{\Phi_p\PAR{\PT{t-s} f}}$  of~\eqref{eq-psi}. A Taylor expansion at $0$ yields
$$
2\rho \, \psi'(0)\leq\psi''(0), 
$$
that is,
$$
 \rho \, \Gamma(h) (x) \leq  \Gamma_2(h) (x) +\frac{2-p}{2(p-1)}\PAR{\frac{\Gamma(h)}{h}}^2 (x)
$$
with $h = (\PT{t} f)^{p-1}$. Then, for instance for $t=0$, $L$ being the Laplacian, so that $\rho=0$, and $f(x) = (x^2+1)^{1/2}$, this writes
$$
0 \leq 1 + \frac{2-p}{2(p-1)} x^4
$$
for all $x$, which leads to a contradiction since $(2-p)/(p-1)<0$.

\medskip

In Theorem \ref{thm-5} below we shall give local $\Phi_p$ entropy inequalities for $p>2$ which are equivalent to the $CD(\rho, \infty)$ criterion.

\bigskip

For $p$ in $]0,1[$, the bound \eqref{eq-phipp} may not hold for any constant $C$ even if the measure $\mu$ satisfies a logarithmic Sobolev inequality. For instance, for the standard Gaussian measure on $\dR$ and $f(x)=e^{\lambda x}$ for $\lambda\in\dR$, \eqref{eq-phipp} writes
$$
e^{\frac{p\lambda^2}{2}(1-p)}-1\leq(1-p)C\lambda^2 \frac{p^2}{4}
$$    
which leads to a contradiction for $\lambda$ going to $+\infty.$

\medskip

In the following theorem we shall also give local $\Phi_p$-entropy inequalities for $p$ in $]0,1[.$

\medskip

\begin{theo}
\label{thm-5}
For $p>0, p \neq 1$ let us consider the function  
$$
\Phi_p(x)=\frac{x^p-x}{p(p-1)}, \quad x >0.
$$ 
 For $\beta \neq 1,2$ we let
\begin{equation}
\label{eq-xi}
\xi(x)=\frac{1-\beta}{2-\beta}\frac{(1+x)^{\frac{2-\beta}{1-\beta}}-1}{x}, \quad x \geq 0
\end{equation}
and $\Delta$ be the set of $(p, \alpha, \beta)$ such that

$\bullet \; \alpha\in ]0,p]$ and $ \beta \in [0,1[ $ if $ p\in ]0,1[, $

$\bullet \; \alpha=1$ and $\beta\geq\frac{4-p}{2-p}$ if $p\in]1,2[, $ 

$\bullet \; \alpha=1$  and $\beta \in \big[ \max \big\{ \frac{p-4}{p-2},0\big\},1 \big[ $ if $p>2.$

\medskip

If $L$ is a diffusion generator and $\rho$ a real number, the following propositions are equiv\nobreak alent, with $\disp \frac{1-e^{-2\rho t}}{2 \rho}$ replaced by $t$ if $\rho=0$:
\begin{enumerate}[(i)]
\item the semigroup $(\Pt)_{t\geq 0}$ satisfies the $CD(\rho,\infty)$ criterion,
\item the semigroup $(\Pt)_{t\geq 0}$ satisfies the local $\Phi_p$-entropy inequality
\begin{equation}
\label{eq-majoration}
\entf{\Pt}^{\Phi_p}\PAR{f}\leq \frac{1-e^{-2\rho t}}{2\rho}\Pt\PAR{f^{p-2}\Gamma(f)}\xi\PAR{\frac{1-e^{-2\rho t}}{2\rho}\kappa_1},
\end{equation}
for all $t\geq 0$, all positive functions $f$ and all $(p,\beta,\alpha)$ in $\Delta$,  where 
$$
\kappa_1={c_p(\beta-1)}\PAR{\frac{\Pt\PAR{f^{\frac{p-b}{\alpha}}\Gamma(f)^{\frac{b}{2\alpha}}}^{\alpha}}{{\Pt\PAR{f^{p-2}\Gamma(f)}}}}^{1-\beta},
$$
$b=2\frac{\beta-2}{\beta-1}$ and $c_p=(2-p)(p-1).$
\item the semigroup $(\Pt)_{t\geq 0}$ satisfies the reverse local $\Phi_p$-entropy inequality
\begin{equation}
\label{eq-minoration}
\entf{\Pt}^{\Phi_p}\PAR{f}\geq\frac{e^{2\rho t}-1}{2\rho} (\Pt f)^{p-2}\Gamma(\Pt f)\xi\PAR{\frac{e^{2\rho t}-1}{2\rho}\kappa_2},
\end{equation}
for all positive functions $f$, all $(p,\beta,\alpha)$ in $\Delta$  and all $t$ in $[0,t_f]$,
where 
$$
\kappa_2={c_p(1-\beta)}e^{-2\rho(2-\beta) t}\PAR{\frac{\Pt(f^{\frac{p-b}{\alpha}}\Gamma(f)^{\frac{b}{2\alpha}})^{\alpha}}{{(\Pt f)^{p-2} \, \Gamma(\Pt f)}}}^{1-\beta},
$$
$b=2\frac{\beta-2}{\beta-1}$, $c_p=(2-p)(p-1)$ and $t_f >0$ depends on $p,\beta,\alpha$ and $f.$
\end{enumerate}
\end{theo}

 Note that for all $(p,\beta,\alpha)\in \Delta$ we have $c_p(1-\beta)\leq0$. 

\begin{ecor}\label{corononadmi} If the semigroup $(\Pt)_{t\geq 0}$  is $\mu$-ergodic and satisfies the  $CD(\rho, \infty)$ criterion with $\rho>0$ then the measure $\mu$ satisfies \begin{equation}
\label{eq-grande2}
\entf{\mu}^{\Phi_p}\PAR{f}\leq \frac{1}{2\rho} \, \mu( f^{p-2}\Gamma(f)) \, \xi\PAR{\frac{1}{2\rho}\kappa_1^\infty},
\end{equation}
for all $(p,\beta,\alpha)\in\Delta$ and positive functions $f$, 
where $\xi$ has been defined in~\eqref{eq-xi} and 
$$
\kappa_1^\infty={c_p(\beta-1)}\SBRA{\frac{\mu(f^{\frac{p-b}{\alpha}}\Gamma(f)^{\frac{b}{2\alpha}})^{\alpha}}{ \mu(f^{p-2}\Gamma(f))}}^{1-\beta}.
$$
\end{ecor}

For admissible functions the bound~\eqref{eq-ad} on the ergodic measure is implied by a weaker averaged criterion as in Propositions~\ref{prop-ci} and~\ref{prop-cip}. This does not seem to be the case for the inequality~\eqref{eq-grande2}, for which our proof is strongly based on the {\it local} $CD(\rho,\infty)$ criterion via the equivalent commutation relation~\eqref{com}.

\begin{erem}
Let us make a few comments on the results of Theorem~\ref{thm-5} and Corollary~\ref{corononadmi} depending on the value of $p$. 
\begin{itemize}
\item If $p\in ]1,2[$, then  all admissible $\beta$ are larger than $2$, so that the map $x\mapsto(1+x)^{\frac{2-\beta}{1-\beta}} $ is strictly concave on $\rr^+$. In particular
$$
\xi(x)=\frac{1-\beta}{2-\beta}\frac{(1+x)^{\frac{2-\beta}{1-\beta}}-1}{x}< 1 
$$
for all $x >0$, which proves that the inequality~\eqref{eq-majoration} is strictly stronger than the local $\Phi_p$-inequality ~\eqref{localPHII} of Section~\ref{sub-11}  with $\Phi = \Phi_p.$ In fact $\disp \lim_{\beta\rightarrow + \infty}\xi(x)=1$ for all $x>0$, so that~\eqref{localPHII} is the limit case of ~\eqref{eq-majoration} as $\beta$ goes to~$+\infty.$ 

Theorem~\ref{thm-main} improved on Theorem~\ref{theophisobolev} by bounding a larger entropy functional by the same energy; here we improve on Theorem~\ref{theophisobolev} by bounding the same entropy functional by a smaller energy.
 Note also that  the method used here does not seem to give the bound~\eqref{eq-ad} of Theorem~\ref{thm-main}. 
 
\item If $p>2$ then all admissible $\beta$ are smaller than $1$, so that  the map  $x\mapsto(1+x)^{\frac{2-\beta}{1-\beta}} $ is convex on $\rr^+$. In particular 
$$
\xi(x)=\frac{1-\beta}{2-\beta}\frac{(1+x)^{\frac{2-\beta}{1-\beta}}-1}{x}\leq (1+x)^{\frac{1}{1-\beta}}
$$
 for all $x>0$, and \eqref{eq-majoration} implies the bound
$$
\entf{\Pt}^{\Phi_p}\PAR{f}\leq\frac{1-e^{-2\rho t}}{2\rho} \Pt\PAR{f^{p-2}\Gamma(f)}\PAR{1+\frac{1-e^{-2\rho t}}{2\rho}\kappa_1}^{\frac{1}{1-\beta}},
$$
where 
$$
\kappa_1={c_p(\beta-1)}\PAR{\frac{\Pt\PAR{f^{{p-b}}\Gamma(f)^{\frac{b}{2}}}}{{\Pt\PAR{f^{p-2}\Gamma(f)}}}}^{1-\beta}.
$$

In particular if  $p\in[2,4]$ then we can let $\beta=0$ to obtain the local inequality
$$
\entf{\Pt}^{\Phi_p}\PAR{f}\leq \frac{1-e^{-2\rho t}}{2\rho} \Pt\PAR{f^{p-2}\Gamma(f)}+\ABS{c_p}\PAR{\frac{1-e^{-2\rho t}}{2\rho}}^2{{\Pt\PAR{f^{{p-4}}\Gamma(f)^{2}}}},
$$
and analogously for the possible ergodic  measure.
\item If $p\in ]0,1[$, then for $\alpha=p$ and $\beta=0$ we obtain the local inequality
\begin{equation}
\label{eq-pp}
\entf{\Pt}^{\Phi_p}\PAR{f}\leq\frac{1-e^{-2\rho t}}{2 \rho} \,  \Pt ( f^{p-2}\Gamma(f)) +\ABS{c_p} \PAR{\frac{1-e^{-2\rho t}}{2\rho}}^2 \Pt( f^{1-\frac{4}{p}}\Gamma(f)^\frac{2}{p})^p
\end{equation}
and analogously for the possible ergodic measure.
\end{itemize}
\end{erem}

\bigskip

Let us now consider the map $\tilde{\Phi}_p$ defined by $\tilde{\Phi}_p(x)=\frac{x^p-1}{p(p-1)}$ for $x >0,$ for which the statements of Theorem~\ref{thm-5} and Corollary~\ref{corononadmi} also hold.  Since $\tilde{\Phi}_p(x)$  tends to $-\log x$ as $p$ goes to $0$, then for instance \eqref{eq-pp} for the ergodic measure leads to the following $\Phi_0$-entropy inequality as $p$ tends to $0$: 
\begin{coro}
 If the semigroup $(\Pt)_{t\geq 0}$  is $\mu$-ergodic and satisfies the  $CD(\rho, \infty)$ criterion with $\rho>0$ then the measure $\mu$ satisfies the inequality
$$
\log\int fd\mu-\int\log f d\mu\leq\frac{1}{2\rho}\int \frac{\Gamma(f)}{f^{2}}d\mu +\frac{1}{2\rho^2}\NRM{ \frac{\Gamma(f)^2}{f^{4}}}_{L^\infty(\mu)}
$$
for all positive functions $f$.
\end{coro}

\bigskip

\noindent
\emph{\textbf{Proof of Theorem~\ref{thm-5}}}.\\\proofbegin~
We first show that $(i)$ implies $(ii)$ and $(iii)$. As in the proof of Theorem~\ref{thm-main} we let
$$
\psi\PAR{s}=\PT{s}\PAR{\Phi_p\PAR{\PT{t-s} f}},
$$
so that
\begin{eqnarray*}
\psi''(s)
&=&
 2 \, \PT{s}\PAR{\frac{\Gamma_2(\Phi_p'(\PT{t-s}   f))}{\Phi''_p(\PT{t-s}   f)}}
+c_p \, \PT{s}\PAR{\PAR{\PT{t-s}   f}^{p-4}{\Gamma (\PT{t-s}   f)}^2} \\
& \geq &
2 \, \rho \, \psi'(s)+c_p \, \PT{s}\PAR{\PAR{\PT{t-s}   f}^{p-4}{\Gamma (\PT{t-s}   f)}^2}
\end{eqnarray*}
by the $CD(\rho,\infty)$ criterion. 

Then the map $(x,y)\mapsto c_p \, x^\beta y^{1-\beta}$ is  convex for all  $(p,\beta,\alpha)$ in $\Delta$, so the second term on the right hand side, which is
$$
c_p\, \PT{s}\PAR{\PAR{\PAR{\PT{t-s} f}^{p-2}{\Gamma (\PT{t-s}   f)}}^{\beta}\PAR{\PAR{\PT{t-s} f}^{p-b}{\Gamma (\PT{t-s}   f)^{b/2}}}^{1-\beta}},
$$
is bounded by below by
$$
c_p\,  \psi'(s)^\beta \, {\PT{s}\PAR{\PAR{\PT{t-s} f}^{p-b}{\Gamma (\PT{t-s}   f)^{b/2}}}}^{1-\beta}.
$$
by the Jensen inequality.

Now, as recalled in Section \ref{subsectPLSI}, the $CD(\rho,\infty)$ criterion is equivalent to the relation 
\begin{equation}\label{com}
\Gamma(\PT{u} f)\leq e^{-2\rho u}\PAR{\PT{u}(\sqrt{\Gamma(f)})}^2
\end{equation}
for all $u \geq 0$ and all positive functions~$f$. Hence
\begin{eqnarray*}
c_p \, {\PT{s}\PAR{\PAR{\PT{t-s} f}^{p-b}{\Gamma (\PT{t-s}   f)^{b/2}}}}^{1-\beta} \! \! \! 
& \geq &
\! \! \! c_p \, e^{-2\rho(2-\beta)(t-s)} {\PT{s}\PAR{\PAR{\PT{t-s} f}^{p-b}\PT{t-s}(\sqrt{\Gamma (f)})^b}}^{1-\beta} \\
& \geq &
\! \! \! c_p \, e^{-2\rho(2-\beta)(t-s)} {\PT{s}\PAR{\PAR{\PT{t-s} f}^{\frac{p-b}{\alpha}}\PAR{\PT{t-s}\sqrt{\Gamma ( {f})}}^{\frac{b}{\alpha}}}}^{\alpha(1-\beta)}
\end{eqnarray*}
by the H\"older inequality.
Now $(x,y)\mapsto x^{\frac{p-b}{\alpha}}y^{\frac{b}{\alpha}}$ is convex for all  $(p,\beta,\alpha)$ in $\Delta$, so that
$$
(\PT{t-s} f)^{\frac{p-b}{\alpha}} \PAR{\PT{t-s}\sqrt{\Gamma ( {f})}}^{\frac{b}{\alpha}}
\leq
 {\PT{t-s}\PAR{f^{\frac{p-b}{\alpha}}{{\Gamma ( {f})}}^{\frac{b}{2\alpha}}}}
$$
by the Jensen inequality, and then
   $$
c_p \, {\PT{s}\PAR{\PAR{\PT{t-s} f}^{\frac{p-b}{\alpha}}\PAR{\PT{t-s}\sqrt{\Gamma ( {f})}}^{\frac{b}{\alpha}}}}^{\alpha(1-\beta)}\geq
c_p \, {\PT{t}\PAR{f^{\frac{p-b}{\alpha}}{{\Gamma ( {f})}}^{\frac{b}{2\alpha}}}}^{\alpha(1-\beta)}
$$
since $c_p(1-\beta)\leq 0$ in all cases. 

Collecting all terms we finally obtain the differential inequality
$$
\psi''(s)\geq2\rho\psi'(s)+c_p e^{-2\rho(2-\beta)(t-s)}\psi'(s)^\beta{\PT{t}\PAR{f^{\frac{p-b}{\alpha}}{{\Gamma ( {f})}}^{\frac{b}{2\alpha}}}}^{\alpha(1-\beta)}, 
 $$
which leads to $(ii)$ by applying the upper bound in Lemma~\ref{lemma-main} below with 
$\mathcal A=c_p {\PT{t}\PAR{f^{\frac{p-b}{\alpha}}{{\Gamma ( {f})}}^{\frac{b}{2\alpha}}}}^{\alpha(1-\beta)}.$
 
 \smallskip
 
Moreover, for $t$ small enough, then  $1+K_2 \geq 0$ in the notation of Lemma~\ref{lemma-main}, and the lower bound in the lemma leads to $(iii)$. 

\bigskip

We now prove that that $(ii)$ implies  $(i)$. For $f = 1 + \varepsilon g$ with $\varepsilon$ going to $0$, then $\kappa_1$ is a $O(\varepsilon^2)$ in the notation of $(ii)$, so that  
$\disp \xi\PAR{\frac{1-e^{-2\rho t}}{2\rho}\kappa_1}$ tends to $1.$ Then $(ii)$ leads to the local Poincar\'e inequality 
$$
\var{\Pt}{g}\leq \frac{1-e^{-2\rho t}}{\rho}\Pt{\Gamma(g)}
$$
in the limit $\varepsilon$ going to $0$, which is equivalent to the $CD(\rho, \infty)$ criterion of $(i)$ as recalled in Section~\ref{subsectPLSI}.

In the same way, $(iii)$ implies the reverse local Poincar\'e inequality 
$$
\var{\Pt}{g}\geq \frac{e^{2\rho t}-1}{\rho}\Gamma(\Pt g )
$$
for all $g$ and $t$, which is also equivalent to the $CD(\rho, \infty)$ criterion of $(i)$; for that purpose we note that $t_{1+\varepsilon g}$ tends  to $+\infty$ as $\varepsilon$ goes to $0$ in the notation of $(iii)$, so that the time limitation in $(iii)$ does not bring any further difficulty.
\proofend

\bigskip

In this proof we have used the following 

\begin{elem}
\label{lemma-main}
Let $\beta\geq0$ with $\beta\neq 1$, $\rho\in\dR$ and $t>0$. Let $\psi$ be a positive, increasing and  $\mathcal C^2$ function on $[0,t]$ such that 
\begin{equation}
\label{eq-lem}
\psi''(s)\geq2\rho\psi'(s)+{\mathcal A} e^{-2\rho(2-\beta)(t-s)}\psi'(s)^\beta
\end{equation}
 for all $s$ in $[0,t]$, where $\mathcal A$ is a real number  such that $\mathcal A(\beta-1)\geq0$. 
Then
\begin{equation}
\label{eq-lem21}
\disp\psi(t)-\psi(0)\leq {\psi'(t)}\frac{1-e^{-2\rho t}}{2\rho}\xi(K_1)
\end{equation}
where $\xi$ is defined in~\eqref{eq-xi} and 
$$
K_1=\frac{1-e^{-2\rho t}}{2\rho}\frac{\mathcal A(\beta-1)}{\psi'(t)^{1-\beta}}. 
$$
If moreover $1+K_2\geq 0$ where 
$$
 K_2=\frac{e^{2\rho t}-1}{2\rho}\frac{\mathcal A(1-\beta)}{\psi'(0)^{1-\beta}}e^{-2\rho(2-\beta)t},
$$
then 
\begin{equation}
\label{eq-lem3}
\psi(t)-\psi(0)\geq {\psi'(0)}\frac{e^{2\rho t}-1}{2\rho}\xi(K_2).
\end{equation}
\end{elem}

\begin{eproof}
We divide equation~\eqref{eq-lem} by $\psi'^{\beta}$, so that
$$
\PAR{\frac{\psi'^{1-\beta}}{1-\beta}}' (s) \geq 2\rho \psi'(s)^{1-\beta}+{\mathcal A} e^{-2\rho(2-\beta)(t-s)}
$$
for all $s$ in $[0,t]$. Then the integration on $[u,v]$ with $0\leq u\leq v \leq t$ gives :
$$
\frac{\psi'(v)^{1-\beta}}{1-\beta}e^{-2\rho(1-\beta)v}-\frac{\psi'(u)^{1-\beta}}{1-\beta}e^{-2\rho(1-\beta)u}\geq
 \frac{\mathcal A}{2\rho }e^{-2\rho(2-\beta)t}(e^{2\rho v}-e^{2\rho u}).
$$
For $u=s$ and $v=t$,  integrating in $s$  over the set $[0,t]$ yields
$$
\psi(t)-\psi(0)\leq \int_0^t e^{-2\rho (t-s)}\PAR{\psi'(t)^{1-\beta}-\frac{\mathcal A (1-\beta)}{2\rho}\PAR{1-e^{-2\rho(t-s)}}}^{\frac{1}{1-\beta}}ds
$$
whether $1- \beta >0$ or $1 - \beta <0,$ which leads to~\eqref{eq-lem21} by the change of variable $x = e^{2 \rho (s-t)}.$

For $u=0$ and $v=s$, we obtain
$$
\frac{\psi'(s)^{1-\beta}}{1-\beta}\geq \frac{e^{2\rho(1-\beta)(s-t)}}{1-\beta} \PAR{\psi'(0)^{1-\beta}e^{2\rho(1-\beta)t}+\frac{\mathcal A (1-\beta)}{2\rho}e^{-2\rho t}\PAR{e^{2\rho s}-1}}
$$
where
$$
\psi'(0)^{1-\beta}e^{2\rho(1-\beta)t}+\frac{\mathcal A (1-\beta)}{2\rho}e^{-2\rho t}\PAR{e^{2\rho s}-1}\geq 0,
$$
for all $s$ in $[0,t]$ since $1 + K_2 \geq0.$ Then integrating over the set $[0,t]$ yields
$$
\psi(t)-\psi(0)\geq \int_0^t e^{-2\rho (t-s)}\PAR{\psi'(0)^{1-\beta}e^{2\rho(1-\beta)t}+\frac{\mathcal A (1-\beta)}{2\rho}e^{-2\rho t}\PAR{e^{2\rho s}-1}}^{\frac{1}{1-\beta}}ds, 
$$
which leads to~\eqref{eq-lem3} by the same change of variable.
\end{eproof}

\section{Applications to a nonlinear evolution equation}
\label{sectappli}

In this section we show how the 
$\Phi$-entropy inequalities studied above for homogeneous Markov semigroups extend to inhomogeneous semigroups and to the solutions of a nonlinear evolution equation in a very simple way.

\medskip

As an example we consider a solution $u = (u_t)_{t \geq 0}$ to the McKean-Vlasov equation
\begin{equation}\label{mkv}
\frac{\partial u_t}{\partial t} = \Delta u_t + \textrm{div} \big( u_t \, \nabla (V + W \ast u_t) \big), \quad t >0, x \in \rr^n
\end{equation}
where $u_t(dx)$ is a probability measure on $\rr^n$ for all $t \geq 0.$ Here $V$ and $W$ are respectively exterior and interaction potentials on $\rr^n$, whereas  $\textrm{div}$ and $\ast$ respectively stand for the divergence and convolution in $x$. This Fokker-Planck type equation has been used in \cite{bccp98} in the modelling of one dimensional granular media in kinetic theory with $V(x) = x^2/2$ and $W(x) = x^3/3$. Explicit rates of convergence  to equilibrium have been obtained in~\cite{cmcv03} in a more general setting.

The convolution term $\nabla W \ast u_t$ induces a nonlinearity in the equation, but nevertheless we shall see how to deduce $\Phi$-entropy inequalities for the solution $u_t$ at time $t$ from those obtained above for diffusion Markov semigroups.

As a first step we derive

\subsection{$\Phi$-entropy inequalities for inhomogeneous Markov semigroups}
\label{sec-inh}

Let $\sigma = (\sigma_{ij})_{1 \leq i,j \leq n}$ be an $n \times n$ matrix and $a(x,y)= (a_i(x,y))_{1 \leq i \leq n}$ have coefficients smooth in $x $ in $\rr^n$ and $y \geq 0.$ For any nonnegative $s$ and $x$ in $\rr^n$ we assume that  the stochastic differential equation
\begin{equation}
\label{eq-last}
dX_t = \sigma \, dB_t - a(X_t, t) \, dt, \quad t \geq s,
\end{equation}
where $(B_t)_{t \geq 0}$ is a standard Brownian motion on $\rr^n,$ has a unique global solution starting from $x$ at time $s$: let it be denoted by $(X_t^{s,x})_{t \geq s}.$ 
One can make the same study for a diffusion matrix $\sigma$ depending on $x$ and $t$, but for simplicity we stick to this simple case, which will be sufficient to be applied to the McKean-Vlasov equation \eqref{mkv}.

Given a function $f$ on $\dR^n$ we let $\PT{s, t} f(x) = \ee f(X_t^{s,x})$ for $0 \leq s \leq t,$ so~that
$$
\frac{\partial}{\partial t} \PT{s, t} f = \PT{s, t} ( L_{(t)} f)  
$$
where  
$$
L_{(t)} f (x) =  \sum_{i,j = 1}^n D_{ij} \,  \frac{\partial^2 f}{\partial x_i \partial x_j} (x) - \sum_{i=1}^n a_i (x,t)  \frac{\partial f}{\partial x_i} (x)
$$
and $D = (D_{ij})_{1 \leq i, j \leq n}$ is the matrix $\disp \frac{1}{2} \sigma \sigma^*.$

To study $\PT{s, t}$ we introduce the following evolution process:  we let  $\bar{\sigma}$ be the $(n+1) \times (n+1)$ matrix with coefficients 
$$
\bar{\sigma}_{ij}= \left\{
\begin{array}{cl}
\sigma_{ij}  & \text{ if } 1 \leq i, j \leq n \\
0 & \text{ otherwise }
\end{array}
\right.
$$
for $1\leq i,j\leq n+1,$ and for $\bar{x}=(x,y)$ in $\rr^n \times \rr^+$ we let $\bar{a} (\bar{x})$ be the vector $(a({x}), 1)$ in~$\rr^{n+1}.$

It follows from our assumptions that for all $(x,y)$ in $\rr^n \times \rr^+$ the stochastic differential equation
$$
 \left\{
\begin{array}{rl}
dX_u =  & \sigma \, dB_u - a(X_u, Y_u) \, du\\
dY_u  = & du
\end{array}
\right. 
$$
with the initial condition $X_0 = x, Y_0 = y$ has a unique global solution on $u \geq 0$, given by $(X_u = X_{y+u}^{y, x}, Y_u = y+u)$ for $u \geq 0$, up to a change of Brownian motion. In other words, for all $\bar{x}$ in $\rr^n \times \rr^+
$ the stochastic differential equation
$$
d\bar{X}_u=   \bar{\sigma} \, d \bar{B}_u - \bar{a}(\bar{X}_u) \, du\
$$
has a unique solution starting from $\bar{x}$ at time $0$: let it be denoted $(\bar{X}_u^{0,\bar{x}})_{u \geq 0}.$ 

Then, for $\bar{f} $ defined on $\rr^n \times \rr^+$ and $u \geq 0$ we let $\bar{\mathbf{P}}_u \bar{f} (\bar{x}) = \ee \bar{f} (\bar{X}_u^{0,\bar{x}}),$ so that  the relation
\begin{equation}\label{barsansbar}
\PT{s, t} f(x) = \bar{\mathbf{P}}_{t-s} \bar{f} (x,s)
\end{equation}
holds for all  functions $f$ on $\rr^n$, $0 \leq s \leq t$ and $x$ in $\rr^n$, where $\bar{f}$ is defined on $\rr^n \times \rr^+$ by $\bar{f} (x,y) = f(x).$

\smallskip
By It\^o's formula, the generator associated to the semigroup $(\bar{\mathbf{P}}_u)_{u \geq 0}$ is given by
$$
\bar{L} \bar{f} (\bar x) =  \sum_{i,j=1}^{n+1} \frac{1}{2} (\bar \sigma \bar \sigma^*)_{ij}  \frac{\partial^2 \bar f}{\partial x_i \partial x_j} (\bar x) - \sum_{i=1}^{n+1} \bar{a}_i (\bar x)  \frac{\partial \bar f}{\partial x_i} ( \bar x).
$$
The diffusion matrix $\bar \sigma \bar \sigma^* /2$ is degenerate, but symmetric and nonnegative, so we are in the setting of the previous sections. In particular, 
according to \eqref{cdrhocst}, the generator $\bar L$ satisfies the $CD(\rho, \infty)$ criterion on $\rr^n \times \rr^+$ if and only if
$$
\frac{1}{2} \big( J \bar a (\bar x) \cdot \frac{1}{2} \bar \sigma \bar \sigma^* + (J \bar a (\bar x) \cdot \frac{1}{2} \bar \sigma \bar \sigma^* )^* \big) \geq \rho \, \frac{1}{2} \bar \sigma \bar \sigma^*
$$
for all $\bar x$ as quadratic forms on $\rr^{n+1}$, that is, if and only if
$$
\frac{1}{2} \big( J a (x,t) \, D + (J a (x,t) \, D )^* \big) \geq \rho \, D
$$
for all $t \geq 0$ and $x$ in $\rr^n$ as quadratic forms on $\rr^n$, hence, if and only if for all $t \geq 0$ the generator $L_{(t)}$ satisfies the $CD(\rho, \infty)$ criterion on $\rr^n$ uniformly on $t$.

Under this assumption, the semigroup $(\bar{\mathbf{P}}_t)_{\geq 0}$ satisfies the $\Phi$-entropy inequalities obtained in Sections \ref{sectphi}, \ref{sect12} and \ref{sectnot12}, for instance for the (simpler) $\Phi$-entropies of Section~\ref{sectphi} with an admissible function $\Phi$ (see the definition in Section~\ref{sub-11}):
 \begin{equation}\label{phibarre}
 \entf{\bar{\PT{u}}}^\Phi (\bar f)\leq \frac{1 - e^{-2\rho u}}{2 \, \rho} \,  \bar{\PT{u}}( \Phi''( \bar f) \bar{\Gamma} (\bar f)) 
 \end{equation}
 for all positive $u$ and all  functions $\bar f = \bar f  (\bar x) $ on $\rr^n \times \rr^+$, where
 $$
 \bar{\Gamma} (\bar f)  = \sum_{i,j = 1}^n \frac{\partial \bar f}{\partial x_i}  D_{ij} \frac{\partial \bar f}{\partial x_j} 
$$
in $\dR^n \!  \times \! \dR^+$ and with $\disp \frac{1 - e^{-2\rho u}}{2 \, \rho}$ replaced by $u$ if $\rho=0.$
Also holds the commutation relation
\begin{equation}\label{combarre}
\bar{\Gamma}(\bar{\PT{u}} \bar f)  \leq e^{-2\rho u} \bar{\PT{u}} \Big(\sqrt{\bar{\Gamma}( \bar f)} \Big)^2 
\end{equation}
 in $\dR^n\times\dR^+$.

\bigskip

Let now $f$ be a given function on $\rr^n$, and let $\bar f$ be defined on $\rr^n \times \rr^+$ by $\bar f (\bar x) = f(x)$ if $\bar x = (x,y)$. Then, by~\eqref{barsansbar},  applying \eqref{phibarre} and \eqref{combarre} to this function $\bar f$ at the point $\bar x = (x,s)$ and $u=t-s$ respectively yield the $\Phi$-entropy inequality
 \begin{equation}\label{phi0t}
 \entf{\PT{s, t}}^\Phi (f) (x) \leq \frac{1 - e^{-2\rho (t-s)}}{2 \, \rho} \,  \PT{s, t}( \Phi''(f) \Gamma (f)) (x), \quad 0 \leq s \leq t, \; x \in \rr^n
 \end{equation}
 and the commutation relation
 \begin{equation}\label{com0t}
\Gamma(\PT{s, t} f) (x) \leq e^{-2\rho (t-s)}{\PT{s, t} \Big(\sqrt{\Gamma(f)} \Big)}^2 (x), \quad 0 \leq s \leq t, \; x \in \rr^n
\end{equation}
for $(\PT{s, t})_{t \geq s \geq 0}$.

Conversely the  $\Phi$-entropy inequality \eqref{phi0t} and the commutation relation \eqref{com0t} independently imply the $CD(\rho, \infty)$ criterion for the generators $L_{(s)}$.

Hence we have obtained:

\begin{prop}\label{inhom}
Let $\Phi$ be an admissible function on an interval $I$. Then, in the above notation, the following three assertions are equivalent, with $\disp  \frac{1 - e^{-2\rho u}}{2 \rho}$ replaced by $u$ if~$\rho =0$:
\begin{enumerate}[(i)]
\item
 the generator $L_{(t)}$ satisfies the ${CD}(\rho,\infty)$ criterion for all nonnegative $t$;
 \item 
 the evolution process $(\PT{s,t})_{0 \leq s \leq t}$ satisfies the commutation relation
 \begin{equation}\label{comst}
\Gamma(\PT{s,t} f) \leq e^{-2\rho (t-s)}{\PT{s,t} \Big(\sqrt{\Gamma(f)} \Big)}^2
\end{equation}
   for all $0 \leq s \leq t$ and all $I$-valued functions $f$ on $\rr^n;$
\item 
 the evolution process $(\PT{s,t})_{0 \leq s \leq t}$ satisfies the local $\Phi$-entropy inequality
 \begin{equation}\label{phist}
 \entf{\PT{s,t}}^\Phi (f) \leq \frac{1 - e^{-2\rho (t-s)}}{2 \, \rho} \,  \PT{s,t}( \Phi''(f) \Gamma (f))
 \end{equation}
    for all $0 \leq s \leq t$ and all $I$-valued functions $f$ on $\rr^n$.
\end{enumerate}
\end{prop}

The implications $(i) \Rightarrow (ii) \Rightarrow (iii)$ have been obtained in \cite{colletmalrieu08} for nonconstant diffusion matrices by rewriting in the inhomogeneous setting the whole argument of Sec\nobreak tion~\ref{sectphi} for homogeneous semigroups. Here the equivalent assertions are obtained without any computation as a simple consequence of the equivalent assertions of Section~\ref{sectphi}, written in a higher dimensional space. 

Let us also note that the $\Phi$-entropy inequalities derived in Sections \ref{sect12} and \ref{sectnot12} can also be transposed to this inhomogeneous setting by the same argument.

\bigskip

In the following section we apply the local bounds of Proposition \ref{inhom} to obtain $\Phi$-entropy inequalities for the solutions of the nonlinear McKean-Vlasov equation \eqref{mkv}.
\subsection{$\Phi$-entropy inequalities for the McKean-Vlasov equation}

In this section we let $(u_t)_{t \geq 0}$ be a solution to~\eqref{mkv} with the probability measure $u_0$ as initial datum.  

Then we let $a(x,t) = \nabla V (x) + \nabla W \ast u_t(x)$ on $\dR^n\times\dR^+$  and we assume that for all $x$ in $\rr^n$ and $s \geq 0$ the stochastic differential inequality~\eqref{eq-last} associated to this $a$ has a unique solution starting from $x$ at time $s$, so that the evolution semigroup $(\PT{s, t})_{t \geq s \geq 0}$ is well defined. 

\begin{theo}\label{theo-mkv}
In the above notation and assumptions, let $\rho$ be a real number, and let $V$ and $W$ be potentials on $\rr^n$ such that $W$ is convex and $\textrm{Hess} \, V (x) \geq \rho \, I$ for all $x$ as quadratic forms on $\rr^n.$ 

If $\Phi$ is an admissible function on an interval  $I$ and $c_0$ a real number, and if the initial datum  $u_0$ satisfies the $\Phi$-entropy inequality with constant $c_0$,
$$
u_0 (\Phi(f)) - \Phi(u_0 (f)) \leq c_0 \, u_0 (\Phi''(f) \Vert \nabla f \Vert^2)
$$
for all $I$-valued functions $f$, then for all $t$ the measure $u_t$ satisfies the $\Phi$-entropy inequality
$$
u_t (\Phi(f)) - \Phi(u_t (f)) \leq \Big(c_0 e^{-2 \rho t} + \frac{1- e^{-2 \rho t}}{2 \rho} \Big) \, u_t (\Phi''(f) \Vert \nabla f \Vert^2)
$$
for all $I$-valued functions $f$, with $\disp  \frac{1 - e^{-2\rho t}}{2 \rho}$ replaced by $t$ if $\rho =0$.
\end{theo}

For $\rho >0$, solutions $u_t$ to~\eqref{mkv} have been shown to converge to a unique equilibrium $u_{\infty}$ as $t$ goes to infinity (see~\cite{cmcv03} for instance). Choosing for instance $u_0$ as a Dirac mass, which satisfies all $\Phi$-entropy inequalities with $c_0 = 0$, and letting $t$ go to infinity in Theorem~\ref{theo-mkv} ensures that for all admissible $\Phi$ the measure $u_{\infty}$ satisfies the $\Phi$-entropy inequality
$$
u_{\infty}(\Phi(f)) - \Phi(u_{\infty}( (f)) \leq \frac{1}{2 \rho} \, u_{\infty}( (\Phi''(f) \Vert \nabla f \Vert^2)
$$
for all maps $f$.

\begin{eproof}
The vector field $a(.,t) = \nabla V + \nabla W \ast u_t$ for $t \geq 0$ is such that
$$
J a(x,t) = \textrm{Hess} \, V (x) + \textrm{Hess} \, W \ast u_t (x) \geq \rho I
$$
for all $x$ and $t$. In particular the generator $L_{(t)} = \Delta - < a (x,t), \nabla >$ satisfies the $CD(\rho, \infty)$ criterion for all $t \geq 0,$ and by $(ii)$ in Proposition \ref{inhom} the local bound
$$
 \PT{0, t}(\Phi (f))  - \Phi( \PT{0, t} f) \leq \frac{1 - e^{-2\rho t}}{2 \, \rho} \,  \PT{0, t}( \Phi''(f) \Vert \nabla f \Vert^2) 
$$
holds for all $f$ and $t$. We now adapt an argument used in \cite{colletmalrieu08} for the propagation of a logarithmic Sobolev inequality by linear evolution equations. 
We integrate with respect to the measure $u_0$ to obtain
$$
u_t(\Phi(f)) = u_0(\PT{0, t}(\Phi (f))) \leq u_0(\Phi( \PT{0, t} f)) + \frac{1 - e^{-2\rho t}}{2 \, \rho} \, u_0 \big(  \PT{0, t}( \Phi''(f) \Vert \nabla f \Vert^2) \big).
$$
On one hand the measure $u_0$ satisfies a $\Phi$-entropy inequality with constant $c_0$ so, letting $g =  \PT{0, t} f$,
$$
u_0(\Phi( \PT{0, t} f)) = u_0(\Phi(g)) \leq \Phi(u_0(g)) + c_0 \, u_0 \big(\Phi''(g) \Vert \nabla g \Vert^2 \big).
$$
First of all $u_0(g) = u_t(f).$ Then $\disp \Vert \nabla g \Vert^2 \leq e^{-2 \rho t} \big( \PT{0, t} \Vert \nabla f \Vert \big)^2$ by $(iii)$ in Proposition \ref{inhom} and the map $(x,y) \mapsto \Phi''(x) \, y^2$ is convex by Remark \ref{remphiadmi}, so by the Jensen inequality
$$
\Phi''(g) \Vert \nabla g \Vert^2 \leq e^{-2 \rho t} \, \Phi''(\PT{0, t} f) \big( \PT{0, t} \Vert \nabla f \Vert \big)^2 \leq e^{-2 \rho t} \, \PT{0, t} \big( \Phi''(f) \Vert \nabla f \Vert^2 \big).
$$
Collecting all terms concludes the argument. 
\end{eproof}

\bigskip

\noindent
{\bf Acknowledgements.} This work was presented at an ``EVOL-Evolution equations" ANR project workshop. It is a pleasure to thank the participants for stimulating discussion on the subject.

\newcommand{\etalchar}[1]{$^{#1}$}


\begin{thebibliography}{AMTU01}

\bibitem[ABC{\etalchar{+}}00]{logsob}
C.~An{\'e}, S.~Blach{\`e}re, D.~Chafa{\"{\i}}, P.~Foug{\`e}res, I.~Gentil,
  F.~Malrieu, C.~Roberto, and G.~Scheffer.
\newblock {\em Sur les in\'egalit\'es de {S}obolev logarithmiques}, volume~10
  of {\em Panoramas et Synth\`eses}.
\newblock Soci\'et\'e Math\'ematique de France, Paris, 2000.

\bibitem[ACJ08]{arnoldcarlenju08}
A.~Arnold, A.~Carlen, and Q.~Ju.
\newblock Large-time behavior of non-symmetric {F}okker-{P}lanck type
  equations.
\newblock {\em Comm. Stoch. Analysis}, 2(1):153--175, 2008.

\bibitem[AD05]{arnolddolbeault05}
A.~Arnold and J.~Dolbeault.
\newblock Refined convex {S}obolev inequalities.
\newblock {\em J. Funct. Anal.}, 225(2):337--351, 2005.

\bibitem[AMTU01]{amtucpde01}
A.~Arnold, P.~Markowich, G.~Toscani, and A.~Unterreiter.
\newblock On convex {S}obolev inequalities and the rate of convergence to
  equilibrium for {F}okker-{P}lanck type equations.
\newblock {\em Comm. Partial Diff. Equations}, 26(1-2):43--100, 2001.

\bibitem[Bak94]{bakrystflour}
D.~Bakry.
\newblock L'hypercontractivit\'e et son utilisation en th\'eorie des
  semigroupes.
\newblock In {\em Lectures on probability theory ({S}aint-{F}lour, 1992)},
  Lecture Notes in Math. 1581, pages 1--114. Springer, Berlin, 1994.

\bibitem[Bak06]{bakrytata}
D.~Bakry.
\newblock Functional inequalities for {M}arkov semigroups.
\newblock In {\em Probability measures on groups: recent directions and
  trends}, pages 91--147. Tata Inst. Fund. Res., Mumbai, 2006.

\bibitem[BCCP98]{bccp98}
D.~Benedetto, E.~Caglioti, J.~A. Carrillo, and M.~Pulvirenti.
\newblock A non-{M}axwellian steady distribution for one-dimensional granular
  media.
\newblock {\em J. Statist. Phys.}, 91(5-6):979--990, 1998.

\bibitem[B{\'E}85]{bakryemery}
D.~Bakry and M.~{\'E}mery.
\newblock Diffusions hypercontractives.
\newblock In {\em S\'eminaire de probabilit\'es, {XIX}, 1983/84}, Lecture Notes
  in Math. 1123, pages 177--206. Springer, Berlin, 1985.

\bibitem[Bec89]{beckner89}
W.~Beckner.
\newblock A generalized {P}oincar\'e inequality for {G}aussian measures.
\newblock {\em Proc. Amer. Math. Soc.}, 105(2):397--400, 1989.

\bibitem[Cha04]{chafai04}
D.~Chafa{\"{\i}}.
\newblock Entropies, convexity, and functional inequalities: on
  {$\Phi$}-entropies and {$\Phi$}-{S}obolev inequalities.
\newblock {\em J. Math. Kyoto Univ.}, 44(2):325--363, 2004.

\bibitem[Cha06]{chafai06}
D.~Chafa{\"{\i}}.
\newblock Binomial-{P}oisson entropic inequalities and the {$M/M/\infty$}
  queue.
\newblock {\em ESAIM Proba. Stat.}, 10:317--339, 2006.

\bibitem[CM08]{colletmalrieu08}
J.-F. Collet and F.~Malrieu.
\newblock Logarithmic {S}obolev inequalities for inhomogeneous semigroups.
\newblock {\em ESAIM Proba. Stat.}, 12:492--504, 2008.

\bibitem[CMV03]{cmcv03}
J.~A. Carrillo, R.~J. McCann, and C.~Villani.
\newblock Kinetic equilibration rates for granular media and related equations:
  entropy dissipation and mass transportation estimates.
\newblock {\em Rev. Mat. Iberoamericana}, 19(3):971--1018, 2003.

\bibitem[DNS08]{dns}
J.~Dolbeault, B.~Nazaret, and G.~Savar\'e.
\newblock {On the Bakry-Emery criterion for linear diffusions and weighted
  porous media equations.}
\newblock {\em Comm. Math. Sci.}, 6(2):477--494, 2008.

\bibitem[EVO09]{evol}
{ANR Project on evolution equations}: EVOL.
\newblock Online encyclopedia available via
  http://www.math.univ-toulouse.fr/evol/categoryencyclopedia, 2009.

\bibitem[Gro75]{gross75}
L.~Gross.
\newblock Logarithmic {S}obolev inequalities.
\newblock {\em Amer. J. Math.}, 97(4):1061--1083, 1975.

\bibitem[Gro93]{gross93}
L.~Gross.
\newblock Logarithmic {S}obolev inequalities and contractivity properties of
  semigroups.
\newblock In {\em Dirichlet forms ({V}arenna, 1992)}, Lecture Notes in Math.
  1563, pages 54--88. Springer, Berlin, 1993.

\bibitem[Hel02]{helffer02}
B.~Helffer.
\newblock {\em Semiclassical analysis, {W}itten {L}aplacians, and statistical
  mechanics}, volume~1 of {\em Series in Partial Differential Equations and
  Applications}.
\newblock World Scientific Publishing Co. Inc., River Edge, 2002.

\bibitem[Led92]{ledoux92}
M.~Ledoux.
\newblock On an integral criterion for hypercontractivity of diffusion
  semigroups and extremal functions.
\newblock {\em J. Funct. Anal.}, 105(2):444--465, 1992.

\bibitem[Led00]{ledouxmarkov}
M.~Ledoux.
\newblock The geometry of {M}arkov diffusion generators.
\newblock {\em Ann. Fac. Sci. Toulouse Math. (6)}, 9(2):305--366, 2000.

\bibitem[LO00]{latalaoles}
R.~Lata{\l}a and K.~Oleszkiewicz.
\newblock Between {S}obolev and {P}oincar\'e.
\newblock In {\em Geometric aspects of functional analysis}, Lecture Notes in
  Math. 1745, pages 147--168. Springer, Berlin, 2000.

\end{thebibliography}

\end{document}